\documentclass[12pt]{article}
\usepackage{fullpage}
\usepackage{amsmath,amsthm}
\usepackage{amssymb,times}
\usepackage{graphicx}
\usepackage{epsfig}
\usepackage{graphicx}
\usepackage{setspace}
\usepackage{natbib}
\newcommand{\silent}[1]{}

\bibpunct{(}{)}{,}{a}{,}{,}
\let\hat\widehat
\let\tilde\widetilde
\newcommand{\inv}[1]{\frac{1}{#1}}
\newcommand{\vecx}{{\bf x}}

\newcommand{\vol}{{\rm Vol}}

\newcommand{\abs}[1]{\left\lvert#1\right\rvert}

\newcommand{\prob}[1]{\ensuremath{\mathbb P}\left(#1\right)}

\newcommand{\norm}[1]{\left\lVert#1\right\rVert}

\newcommand{\beq}{\begin{equation}}
\newcommand{\eeq}{\end{equation}}
\newcommand{\ben}{\begin{eqnarray*}}
\newcommand{\een}{\end{eqnarray*}}

\newcommand{\X}{{\mathcal X}}

\newtheorem{theorem}{Theorem}[section]

\newtheorem{lemma}[theorem]{Lemma}

\newtheorem{definition}[theorem]{Definition}
\newtheorem{example}[theorem]{Example}

\newtheorem{remark}[theorem]{Remark}

\def\qed{\hskip1pt $\;\;\scriptstyle\Box$}
\newenvironment{proofof}[1]{\hspace*{20pt}{\it Proof}{ of #1}.\hskip10pt}{\qed\vskip5pt}
\newenvironment{proofof2}{\hskip10pt}{\qed\vskip5pt}

\def\mld{${}^\ddag$}
\def\dos{${}^*$}
\def\sfs{${}^\dag$}
\begin{document}

\begin{center}
{\bf\Large A Statistical Framework for Differential Privacy}\footnote{
We thank Avrim Blum, Katrina Ligett,
Steve Fienberg, Alessandro Rinaldo and Yuval Nardi
for many helpful discussions.
We thank Wenbo Li and Mikhail Lifshits for helpful pointers and discussions 
on small ball probabilities.
We thank the Associate Editor and three referees for a plethora of comments
that led to improvements in the paper.
Research supported by NSF grant CCF-0625879, a Google research grant and a grant from Carnegie Mellon's Cylab. 
The second author is also partially supported by the Swiss National Science Foundation (SNF) Grant 20PA21-120050/1.
}\\
\vskip .2in
\begin{tabular}{c}
{\large
 Larry Wasserman\dos\mld\;\; Shuheng Zhou\sfs} \\[15pt]
{\dos}Department of Statistics \\
{\mld}Machine Learning Department \\
Carnegie Mellon University \\
Pittsburgh, PA 15213 \\[20pt]
{\sfs}Seminar f\"{u}r Statistik \\
ETH Z\"{u}rich, CH 8092 \\[20pt]
\today \\[5pt]
\end{tabular}
\end{center}

\begin{quote}
  One goal of statistical privacy research is to construct a data
  release mechanism that protects individual privacy while preserving
  information content.  An example is a {\em random mechanism} that
  takes an input database $X$ and outputs a random database $Z$
  according to a distribution $Q_n(\cdot |X)$.  {\em Differential
  privacy} is a particular privacy requirement developed by computer
  scientists in which $Q_n(\cdot |X)$ is required to be insensitive to
  changes in one data point in $X$.  This makes it difficult to infer
  from $Z$ whether a given individual is in the original database $X$.
  We consider differential privacy from a statistical perspective.  We
  consider several data release mechanisms that satisfy the
  differential privacy requirement.  We show that it is useful to
  compare these schemes by computing the rate of convergence of
  distributions and densities constructed from the released data.  We
  study a general privacy method, called the exponential mechanism,
  introduced by~\cite{MT07}.  We show that the accuracy of this method
  is intimately linked to the rate at which the probability that the
  empirical distribution concentrates in a small ball around the true
  distribution.
\end{quote}

\setstretch{1.5}

\section{Introduction}

One goal of data privacy research is to derive a mechanism that takes
an input database $X$ and releases a transformed database $Z$ such
that individual privacy is protected yet information content is
preserved.  This is known as disclosure limitation.  In this paper we
will consider various methods for producing a transformed database $Z$
and we will study the accuracy of inferences from $Z$ under various
loss functions.

There are numerous approaches to this problem.  The literature is vast
and includes papers from computer science, statistics and other
fields.  The terminology also varies considerably. We will use the
terms ``disclosure limitation'' and ``privacy guarantee''
interchangeably.

Disclosure limitation methods include clustering
\citep{sweeney2002kam,AFK+06}, $\ell$-diversity \citep{kifer},
$t$-closeness \citep{ninghui}, data swapping \citep{fienberg2004dsv},
matrix masking \citep{Ting:07}, cryptographic approaches
\citep{pinkas2002ctp, feigen:06}, data
perturbation~\citep{evfimievski2004ppm,kim2003mnm,warner1965rrs,Fienberg:98}
and distributed database methods~\citep{Fienberg:07,Sanil:04}.
Statistical references on disclosure risk and limitation include
\cite{DuncanLambert:86,DuncanLambert:89,duncan:91,Reiter:05}.  We
refer to~\cite{Reiter:05} and~\cite{Sanil:04} for further references.

One approach to defining a privacy guarantee that has received much
attention in the computer science literature is known as {\em
differential privacy}~\citep{DMNS06,Dwork:06}.  There is a large
body of work on this topic including, for example,
\cite{DN03,DN04,BDMN05,DMT07,NRS07,barak2007paa,MT07,BLR08,KLN+08}.
\cite{BLR08} gives a machine learning approach to
inference under differential privacy constraints and
to some extent
our results are inspired by that paper.
\cite{Smith08} shows how to provide efficient point estimators
while preserving differential privacy.
He constructs estimators for parametric models
with mean squared error
$(1+o(1))/(nI(\theta))$ where
$I(\theta)$ is the Fisher information.
\cite{MKAGV08} consider privacy for histograms
by sampling from the posterior distribution of the
cell probabilities.
We discuss \cite{MKAGV08} further in Section \ref{sec::histogram}.
After submitting the first draft of this paper,
new work has appeared on differential privacy
that is also statistical in nature, namely,
~\cite{GRS09,DL09,DNRRV09,DFKN09}.

The goals of this paper
are to explain differential privacy
in statistical language,
to show how to compare
different privacy mechanisms by
computing the rate of convergence of
distributions and densities based on the released data $Z$,
and to study a general privacy method, called the exponential mechanism,
due to~\cite{MT07}.
We show that the accuracy of this method is intimately linked to
the rate at which the probability that the empirical distribution
concentrates in a small ball around the true distribution.
These so called ``small ball probabilities''
are well-studied in probability theory.
To the best of our knowledge, this
is the first time a connection has been made between differential privacy and
small ball probabilities.
We need to make two disclaimers.
First, the goal of our paper is to investigate differential privacy.
We will not attempt to review all approaches to privacy
or to compare differential privacy with other approaches.
Such an undertaking is beyond the scope of this paper.
Second, we focus only on statistical properties here.
We shall not concern ourselves in this paper with computational efficiency.

In Section \ref{sec::diffpriv}
we define differential privacy and provide motivation
for the definition.
In Section \ref{sec::informative}
we discuss conditions that
ensure that a privacy mechanism
preserves information.
In Section \ref{sec::histogram}
we consider two histogram based methods.
In Section \ref{sec::exponential} and \ref{sec::density},
we examine another method known as the
exponential mechanism.
Section \ref{sec::simulations} contains
a small simulation study and
Section \ref{sec::conclusion} contains concluding remarks.
All technical proofs appear in Section~\ref{sec::proofs}. 

\subsection{Summary of Results}

We consider several different data release mechanisms that satisfy
differential privacy.
We evaluate the utility of these mechanisms
by evaluating the rate at which $d(P,P_Z)$ goes to 0,
where $P$ is the distribution of the data $X\in {\cal X}$,
$P_Z$ is the empirical distribution of the released data $Z$,
and $d$ is some distance between distributions.
This gives an informative way to compare data release mechanisms.
In more detail, we consider the Kolmogorov-Smirnov (KS) distance:
$\sup_{x \in \X} |F(x) - \hat{F}_Z(x)|$,
where $F$, $\hat{F}_Z$ denote the cumulative distribution 
function (cdf) corresponding to $P$ and the empirical 
distribution function corresponding to $P_Z$, respectively.
We also consider the squared $L_2$ distance: 
$\int (p(x) - \hat{p}_Z)^2$, where $\hat{p}_Z$ is a density 
estimator based on $Z$.
Our results are summarized in the following tables, where $n$ 
denotes the sample size.

The first table concerns the case where the data are in $\mathbb{R}^r$
and the density $p$ of $P$ is Lipschitz.
Also reported are the minimax rates  of convergence for density 
estimators in KS and in squared $L_2$ distances.
We see that the accuracy depends both on the data releasing mechanism
and the distance function $d$.
The results are from Sections 4 and 5 of the paper.
(The exponential mechanism under $L_2$ distance is marked NA
but is in the second table in case $r=1$. We note that
the rate for KS distance for perturbed histogram 
is $\sqrt{\log n/n}$ for $r=1$.)

\begin{center}
\begin{tabular}{|l||l|l|l||l|}\hline
                      & \multicolumn{3}{c||} {Data Release mechanism} & \\
\cline{2-4}
Distance    & smoothed    & perturbed   & exponential & minimax \\
            & histogram   & histogram   & mechanism   & rate \\ \hline
$L_2$       & $n^{-2/(2r+3)}$ & $n^{-2/(2+r)}$  & NA    & $n^{-2/(2+r)}$ \\
Kolmogorov-Smirnov &  $\sqrt{\log n} \times n^{-2/(6+r)}$  & 
$\log n \times n^{-2/(2+r)}$
& $n^{-1/3}$        & $n^{-1/2}$\\ \hline
\end{tabular}
\end{center}

The next table summarizes the results for 
the case where the dimension of $X$ is $r=1$ and 
the density $p$ is assumed to be in a Sobolev space of 
order $\gamma$. 
We only consider the squared $L_2$ distance between the true density 
$p$ and the estimated density $\hat{p}_Z$ in this case.
The results are from Section 6 of the paper.

\begin{center}
\begin{tabular}{|l||l|l||l|}\hline
      & exponential   & perturbed orthogonal & minimax rate\\ 
      & mechanism     & series estimator & \\ \hline
$L_2$ & $n^{-\gamma/(2\gamma+1)}$ & $n^{-2\gamma/(2\gamma+1)}$ & $n^{-2\gamma/(2\gamma+1)}$\\ \hline
\end{tabular}
\end{center}

Our results show that, in general, privacy schemes seem not to
yield minimax rates.
Two exceptions are perturbation methods evaluated under $L_2$ loss
which do yield minimax rates.
An open question is whether the slower than minimax rates
are intrinsic to the privacy methods.
It is possible, for example, that our rates are not tight.
This question could be answered by establishing lower bounds
on these rates. We consider this an important topic for future research.

\section{Differential Privacy}
\label{sec::diffpriv}

Let $X_1,\ldots, X_n$ be a random sample 
(independent and identically distributed) of size $n$ from a distribution $P$
where $X_i \in {\cal X}$.
To be concrete, we shall assume that
${\cal X}\equiv [0,1]^r = [0,1]\times [0,1] \times \cdots \times [0,1]$
for some integer $r \geq 1$.
Extensions to more general sample spaces are certainly possible
but we focus on this sample space to avoid unnecessary technicalities.
(In particular, it is difficult to extend differential privacy to unbounded domains.)
Let $\mu$ denote Lebesgue measure and let
$p=dP/d\mu$ if the density exists.
We call $X= (X_1,\ldots, X_n)$ a database.
Note that $X\in {\cal X}^n = [0,1]^r \times \cdots \times [0,1]^r$.
We focus on mechanisms
that take a database $X$ as input
and output a sanitized database 
$Z = (Z_1,\ldots, Z_k)\in {\cal X}^k$
for public release.
In general, $Z$ need not be the same size as $X$.
For some schemes,
we shall see that large $k$ can lead to low privacy and high accuracy while
while small $k$ can lead to high privacy and low accuracy.
We will let $k \equiv k(n)$ change with $n$.
Hence, any asymptotic statements involving $n$ increasing will
also allow $k$ to change as well.

A {\em data release mechanism}
$Q_n(\cdot|X)$ is a conditional distribution 
for $Z=(Z_1,\ldots, Z_k)$ given $X$.
Thus,
$Q_n(B|X=x)$
is the probability that the output database
$Z$ is in a set $B\in {\cal B}$ given that the input database is $x$,
where ${\cal B}$ are the measurable subsets of ${\cal X}^k$.
We call
$Z=(Z_1,\ldots, Z_k)$
a {\em sanitized database}.
Schematically:
$$
{\rm input\ database\ }X=(X_1,\ldots, X_n)\ 
\xrightarrow[\rm sanitize]{Q_n(Z|X)}
{\rm output\ database\ }Z=(Z_1,\ldots, Z_k).
$$
The marginal distribution of
the output database $Z$ induced by $P$ and $Q_n$ is
$M_n(B) = \int Q_n(B|X=x)dP^n(x)$
where $P^n$ is the $n$-fold product measure of $P$.

\begin{example}
A simple example to help the reader have a concrete example in mind
is adding noise.
In this case, $Z = (Z_1,\ldots, Z_n)$
where $Z_i = X_i + \epsilon_i$
and
$\epsilon_1,\ldots,\epsilon_n$ are mean 0 independent
observations drawn from
some known distribution $H$ with density $h$.
Hence
$Q_n$ has density
$q_n(z_1,\ldots, z_n|x_1,\ldots, x_n) = \prod_{i=1}^n h(z_i - x_i)$.
\end{example}

\begin{definition}
Given two databases $X=(X_1,\ldots, X_n)$ and
$Y=(Y_1,\ldots, Y_n)$,
let $\delta(X,Y)$ denote the Hamming distance between 
$X$ and $Y$: $\delta(X,Y) = \#\Bigl\{ i:\ X_i \neq Y_i \Bigr\}$.
\end{definition}

A general data release mechanism 
is the {\em exponential mechanism}~\citep{MT07} which is defined as follows.
Let $\xi: {\cal X}^n \times {\cal X}^k:\to [0,\infty)$ be any function.
Each such $\xi$ defines a different exponential mechanism.
Let
\begin{equation}
\Delta \equiv \Delta_{n,k} = 
\sup_{\stackrel{x,y\in {\cal X}^n}{\delta(x,y)=1}} 
\sup_{z\in {\cal X}^k}
|\xi(x,z) - \xi(y,z)|,
\end{equation}
that is, $\Delta_{n,k}$ is the maximum change to $\xi$ 
caused by altering a single entry in $x$.
Finally, let
$(Z_1,\ldots, Z_k)$ 
be a random vector drawn from
the density
\begin{equation}
h(z|x) = 
\frac{\exp\left( -\frac{\alpha \xi(x,z)}{2\Delta_{n,k}}\right)}
     {\int_{{\cal X}^k}\exp\left( -\frac{\alpha \xi(x,s)}{2\Delta_{n,k}}\right)ds}
\end{equation}
where $\alpha \geq 0$,
$z = (z_1,\ldots, z_k)$ and
$x=(x_1,\ldots,x_n)$.
In this case, $Q_n$ has density $h(z|x)$.
We'll discuss the exponential mechanism 
in more detail later.

There are many definitions of privacy
but in this paper we focus on the
following 
definition due to 
\cite{DMNS06} and
\cite{Dwork:06}.

\begin{definition}
Let $\alpha \geq 0$.
We say that
$Q_n$ satisfies
{\em $\alpha$-differential privacy}
if
\begin{equation}
\label{eq::diffpriv}
\sup_{\stackrel{x,y\in {\cal X}^n}{ \delta(x,y)=1}}
\sup_{B\in {\cal B}} 
\frac{Q_n(B|X=x)}{Q_n(B|X=y)}\leq e^\alpha
\end{equation}
where ${\cal B}$ are the measurable sets on ${\cal X}^k$.
The ratio is interpreted to be 1 whenever
the numerator and denominator are both 0.
\end{definition}

The definition of differential privacy is based on ratios of probabilities.
It is crucial to measure closeness by ratios of probabilities since that 
protects rare cases which have small probability under $Q_n$.
In particular, if changing one entry in the database $X$ cannot change the 
probability distribution $Q_n(\cdot|X=x)$ very much, then we can claim that
a single individual cannot guess whether he is in the original database 
or not. The closer $e^\alpha$ is to 1, the stronger privacy 
guarantee is. Thus, one typically chooses $\alpha$ close to 0.
See~\cite{DMNS06} for more discussion on these points.
Indeed, suppose that two subjects each believe that one of them is in the 
original database.
Given $Z$ and full knowledge of 
$P$ and $Q_n$ can they test who is in $X$?
The answer is given in the following result.
(In this result, we drop the assumption that the user does not know $Q_n$.)

\begin{theorem}\label{thm::power}
Suppose that $Z$ is obtained from a data release mechanism
that satisfies $\alpha$-differential privacy.
Any level $\gamma$ test
which is a function of $Z$, $P$ and $Q_n$
of $H_0: X_i = s$ versus 
$H_1: X_i = t$
has power bounded above by
$\gamma e^{\alpha}$.
\end{theorem}

Thus, if $Q_n$ satisfies differential privacy
then it is virtually impossible to test the hypothesis
that either of the two subjects is in the database
since the power of such a test is nearly equal to its level.
A similar calculation shows that if one does a Bayes test between
$H_0$ and $H_1$ then the Bayes factor
is always between
$e^{-2\alpha}$ and $e^{2\alpha}$.
For more detail on the motivation for the definition as well as consequences,
see \cite{DMNS06,Dwork:06,GKS08,Ras09}.

The following result,
which is proved in \cite{MT07} (Theorem~6),
shows that the exponential mechanism always preserves differential privacy.

\begin{theorem}\textnormal{\citep{MT07}}
The exponential mechanism satisfies the $\alpha$-differential privacy.
\end{theorem}

To conclude this section we record a few useful facts.
Let $T(X,R)$ be a function of $X$ and some auxiliary random 
variable $R$ which is independent of $X$.
After including this auxiliary random variable
we define differential privacy as before.
Specifically, $T(X,R)$ satisfies differential privacy if for all $B$,
and all $x,x'$ with $\delta(x,x')=1$ we have that
$\mathbb{P}(T(X,R)\in B |X=x) \leq e^\alpha \mathbb{P}(T(X,R)\in B |X=x')$.
The third part is Proposition 1 from~\cite{DMNS06}.

\begin{lemma}\label{lemma::useful}
We have the following:
\begin{enumerate}
\item If $T(X,R)$ satisfies differential privacy 
then $U = h(T(X,R))$ also satisfies
differential privacy for any measurable function $h$.
\item Suppose that $g$ is a density function
constructed from a random vector $T(X,R)$ that 
satisfies differential privacy.
Let $Z=(Z_1, \ldots, Z_k)$ be $k$ iid draws from $g$.
This defines a mechanism
$Q_n(B|X) = \mathbb{P}(Z\in B|X)$.
Then $Q_n$ satisfies differential privacy
for any $k$.
\item (Proposition 1 from \cite{DMNS06}.)
Let $f(x)$ be a function of $x=(x_1, \ldots, x_n)$
and define
$S(f) = \sup_{x,x': \delta(x,x')=1} \norm{f(x) - f(x')}_1$
where $\norm{a}_1 = \sum_j |a_j|$.
Let $R$ have density
$g(r) \propto e^{-\alpha |r|/S(f)}$.
Then $T(X,R) = f(X) + R$ satisfies differential privacy.
\end{enumerate}
\end{lemma}

\section{Informative Mechanisms}
\label{sec::informative}

A challenge in privacy theory
is to find
$Q_n$ that satisfies differential privacy and yet
yields datasets $Z$ that preserve information.
Informally, a mechanism is informative
if it is possible to make precise inferences from the released data
$Z_1, \ldots, Z_k$.
Whether or not a mechanism is informative will depend on the goals of the 
inference. From a statistical perspective,
we would like to infer $P$ or functionals of $P$
from $Z$. 
\cite{BLR08}
show that
the probability content of some classes
of intervals can be estimated accurately
while preserving privacy.
Their results motivated the current paper.
We will assume throughout that the user has access to the
sanitized data $Z$ but not the mechanism $Q_n$.
The question of how a data analyst can use knowledge
of $Q_n$ to improve inferences is left to future work.

There are many ways to measure the information in $Z$.
One way is through distribution functions.
Let $F$ denote the cumulative distribution function (cdf)
on ${\cal X}$ corresponding to $P$.
Thus
$F(x) = P( X\in (-\infty,x_1]\times \cdots \times (-\infty,x_r] )$
where $x=(x_1,\ldots, x_r)$.
Let $\hat{F}\equiv\hat{F}_X$ denote the empirical distribution function
corresponding to $X$ and similarly let
$\hat{F}_Z$ denote the empirical distribution function
corresponding to $Z$.
Let $\rho$ denote any distance measure on distribution functions.

\begin{definition}
$Q_n$ is consistent with respect to $\rho$ if
$\rho(F,\hat{F}_Z)\stackrel{P}{\to} 0$.
$Q_n$ is $\epsilon_n$-informative if
$\rho(F,\hat{F}_Z) = O_P(\epsilon_n)$.
\end{definition}

\vspace{.5cm}

An alternative to requiring $\rho(F,\hat{F}_Z)$ to be small
is to require
$\rho(\hat{F},\hat{F}_Z)$ to be small.
Or one could require 
$Q_n( \rho(\hat{F},\hat{F}_Z) > \epsilon |X=x)$ be small
for all $x$
as in \cite{BLR08}.
These requirements are similar.
Indeed, suppose $\rho$ satisfies the triangle inequality
and that
$\hat{F}$ is consistent in the $\rho$ distance, that is,
$\rho(\hat{F},F)\stackrel{P}{\to}0$.
Assume further that $\rho(\hat{F}, F) = O_P(\epsilon_n)$.
Then 
$\rho(F,\hat{F}_Z) = O_P(\epsilon_n)$
implies that
$$
\rho(\hat{F},\hat{F}_Z) \leq \rho(\hat{F},F) + \rho(F,\hat{F}_Z) =
O_P(\epsilon_n);
$$
Similarly,
$\rho(\hat{F},\hat{F}_Z) = O_P(\epsilon_n)$
implies that
$\rho(F,\hat{F}_Z) = O_P(\epsilon_n)$.

Let $\mathbb{E}_{P,Q_n}$
denote the expectation under the joint distribution
defined by $P^n$ and $Q_n$.
Sometimes we write
$\mathbb{E}$ when there is no ambiguity.
Similarly,
we use $\mathbb{P}$ to denote
the marginal probability under $P^n$ and $Q_n$:
$\mathbb{P}(A) = \int_A dQ_n(z_1,\ldots, z_k|x_1,\ldots, x_n) dP(x_1)\cdots dP(x_n)$ for $A\in {\cal X}^k$.

There are many possible choices for $\rho$.
We shall mainly focus on
the Kolmogorov-Smirnov (KS) distance
$\rho(F,G)= \sup_x |F(x) - G(x)|$
and the squared $L_2$ distance
$\rho(F,G)= \int (f(x) - g(x))^2 dx$
where 
$f= dF/d\mu$ and
$g= dG/d\mu$.
However, our results can be carried over to 
other distances as well.

Before proceeding let us note that
we will need some assumptions on $F$
otherwise we cannot have a consistent scheme
as shown in the following theorem.
The following result --- essentially a 
re-expression of a result in \cite{BLR08}
in our framework --- makes this clear.

\begin{theorem}\label{thm::nogo}
Suppose that $Q_n$ satisfies differential privacy and
that $\rho(F,G) = \sup_x |F(x) - G(x)|$.
Let $F$ be a point mass distribution.
Thus $F(y) = I(y \geq x)$ for some point $x\in [0,1]$.
Then $\hat{F}_Z$ is inconsistent, that is,
there is a $\delta>0$ such that
$\liminf_{n\to\infty}P^n(\rho(F,\hat{F}_Z)> \delta) > 0$.
\end{theorem}

\section{Sampling From a Histogram}
\label{sec::histogram}

The goal of this section is to give two concrete, simple 
data release methods that achieve differential privacy.
The idea is to draw a random sample from histogram.
The first scheme draws observations from a smoothed histogram.
The second scheme draws observations from a randomly perturbed histogram.
We use the histogram for its familiarity and simplicity and because it is
used in applications of differential privacy.
We will see that the histogram has to be carefully constructed
to ensure differential privacy.
We then compare the two schemes
by studying the accuracy of the inferences from the released data.
We will see that the accuracy depends both on
how the histogram is constructed and on what measure of
accuracy we use.

Let $L>0$ be a constant and
suppose that $p=dP/d\mu \in {\cal P}$
where
\begin{equation}
\label{eq::lipschitz}
{\cal P} = \Biggl\{ p:\ |p(x) - p(y)| \leq L ||x-y|| \Biggr\}
\end{equation}
is the class of Lipschitz functions.
We assume throughout this section that $p\in {\cal P}$.
The minimax rate of convergence for density estimators in squared 
$L_2$ distance for ${\cal P}$ is $n^{-2/(2+r)}$~\citep{scot:1992}.

Let $h=h_n$
be a binwidth such that
$0 < h < 1$ and such that $m = 1/h^r$ is an integer.
Partition ${\cal X}$ into
$m$ bins 
$\{B_1,\ldots, B_m\}$
where each bin $B_j$ is a cube with sides of length $h$.
Let $I(\cdot )$ denote the indicator function.
Let $\hat{f}_m$ denote the corresponding histogram estimator
on ${\cal X}$, namely,
$$
\hat{f}_{m}(x) = \sum_{j=1}^m \frac{\hat{p}_j}{h^r}I(x\in B_j)
$$
where
$\hat{p}_j = C_j/n$
and
$C_j = \sum_{i=1}^n I(X_i\in B_j)$
is the number of observations in $B_j$.
Recall that $\hat{f}_m$ is a consistent estimator of $p$
if $h = h_n \to 0$ and
$n h_n^r \to \infty$.
Also, the optimal choice of
$m=m_n$ for $L_2$ error
under ${\cal P}$
is $m_n \asymp n^{r/(2+r)}$, in which case $\int (p-\hat{f}_m)^2 = O_P(n^{-2/(2+r)})$~\citep{scot:1992}.
Here, $a_n\asymp b_n$ means that both
$a_n/b_n$ and $b_n/a_n$ are bounded for large $n$.

\subsection{Sampling from a Smoothed Histogram}
The first method for generating released data $Z$ from a histogram while
achieving differential privacy proceeds as follows. 
Recall that the sample space is $[0,1]^r$.
Fix a constant $0 < \delta < 1$ and define the smoothed histogram
\begin{equation}
\label{eq::reasoning}
\hat{f}_{m,\delta}(x) = (1-\delta)\hat{f}_{m}(x) + \delta.
\end{equation}

\begin{theorem}
\label{thm::hist-theorem}
Let
$Z = (Z_1,\ldots, Z_k)$ where
$Z_1,\ldots, Z_k$ are $k$ iid draws from
$\hat{f}_{m,\delta}(x)$.
If
\begin{equation}
\label{eq::needed}
k \log \left( \frac{(1-\delta)m}{n\delta} + 1\right) \leq \alpha
\end{equation}
then $\alpha$-differential privacy holds.
\end{theorem}

Note that for $\delta \to 0$ and $\frac{m}{n\delta} \to 0$,
$\log \left( \frac{(1-\delta)m}{n\delta} + 1\right) =
\frac{m}{n\delta} (1 + o(1)) \approx \frac{m}{n\delta}$.
Thus~\eqref{eq::needed} is approximately the same as requiring
\begin{equation}
\label{eq::smple}
\frac{mk}{\delta}  \leq  n \alpha.
\end{equation}


Equation (\ref{eq::smple})
shows an interesting tradeoff between $m$, $k$ and $\delta$.
We note that
sampling from the usual
histogram corresponding to $\delta=0$ does not preserve differential privacy.

Now we consider how to choose $m,k,\delta$ to minimize
$\mathbb{E}(\rho(F,\hat{F}_Z))$ while
satisfying \eqref{eq::needed}.
Here, $\mathbb{E}$ is the expectation under the randomness
due to sampling from $P$ and due to the privacy mechanism
$Q_n$.
Thus, for any measurable function $h$,
$$
\mathbb{E}(h(Z)) =
\int\int h(z_1,\ldots, z_k) dQ_n(z_1,\ldots, z_k|x_1,\ldots, x_n) dP(x_1)\cdots dP(x_n).
$$

Now we give a result that shows how accurate the inferences are
in the KS distance using the smoothed histogram sampling scheme.

\begin{theorem}\label{thm::histogram}
Suppose that $Z_1,\ldots, Z_k$ are drawn as described in the previous theorem.
Suppose~\eqref{eq::lipschitz} holds.
Let $\rho$ be the KS distance.
Then choosing
$m \asymp n^{r/(6+r)}$, $k \asymp m^{4/r}=n^{ 4/(6+r)}$ and $\delta = (mk/n\alpha)$
minimizes $\mathbb{E}\rho(F,\hat{F}_Z)$ 
subject to \eqref{eq::needed}.
In this case, $\mathbb{E}\rho(F,\hat{F}_Z) = O\left(\frac{\sqrt{\log n}}{n^{2/(6+r)}}\right).$
\end{theorem}

In this case we see that we have consistency since
$\rho(F,\hat{F}_Z)=o_P(1)$ but the rate is slower than the minimax rate of 
convergence for density estimators in KS distance, which is $n^{-1/2}$. 
Now let $\hat{q}_j = \#\{Z_i \in B_j\}/k$ and
\begin{equation}\label{eq::lwidtsdefn}
\rho(F,\hat{F}_Z) = \int (p(x) - \hat{f}_Z(x))^2 dx, 
\; \text{where } \; 
\hat{f}_Z(x) = h^{-r}\sum_{j=1}^m \hat{q}_j I(x\in B_j).
\end{equation}

\begin{theorem}\label{thm::histogram2}
Assume the conditions of the previous theorem. Let $\rho$ be the 
squared $L_2$ distance as defined in (\ref{eq::lwidtsdefn}).
Then choosing
\begin{equation*}
m \asymp n^{r/(2r+3)},\ \ \ 
k \asymp n^{(r+2)/(2r+3)},\ \ \ 
\delta \asymp n^{-1/(r+3)}
\end{equation*}
minimizes $\mathbb{E}\rho(F,\hat{F}_Z)$ 
subject to \eqref{eq::needed}.
In this case, $\mathbb{E}\rho(F,\hat{F}_Z) = O(n^{-2/(2r+3)})$.
\end{theorem}
Again, we have consistency
but the rate is slower than the minimax rate which is $n^{-2/(2+r)}$.
\citep{scot:1992}

\subsection{Sampling From a Perturbed Histogram}
The second method, which we call the
sampling from a perturbed histogram, is due to
Dwork et. al. (2006).
Recall that $C_j$ is the number of observations in bin $B_j$.
Let 
$D_j = C_j + \nu_j$
where
$\nu_1, \ldots, \nu_m$ are independent, identically
distributed draws from a Laplace density
with mean 0 and variance $8/\alpha^2$.
Thus the density of $\nu_j$ is
$g(\nu) = (\alpha/4)e^{-|\nu|\alpha/2}$.
Dwork et. al. (2006) show that
releasing $D=(D_1, \ldots, D_m)$
preserves differential privacy.
However, our goal is to release a database $Z=(Z_1,\ldots, Z_k)$
rather than just a set of counts.
Now define
$$\tilde{D}_j = \max\{ D_j,0\} \; \text{ and } \; 
\hat{q}_j = \tilde{D}_j/\sum_s \tilde{D}_s.$$
Since $D$ preserves differential privacy, it follows 
from Lemma \ref{lemma::useful} that
$(\hat{q}_1, \ldots, \hat{q}_m)$ also
preserve differential privacy;
Moreover, any sample  $Z = (Z_1, \ldots, Z_k)$
from $\tilde{f}(x) = h^{-r}\sum_{j=1}^m \hat{q}_jI(x\in B_j)$
preserve differential privacy
for any $k$.

\begin{theorem}\label{thm::histogram3}
Let $Z = (Z_1, \ldots, Z_k)$
be drawn from $\tilde{f}(x) = h^{-r}\sum_{j=1}^m \hat{q}_jI(x\in B_j)$.
Assume that there exists a constant $1 \leq C < \infty$ such that 
$\sup_x p(x) = C$.\\
(1) Let $\rho$ be the $L_2$ distance and
$\hat{f}_Z$ be as defined in (\ref{eq::lwidtsdefn}).
Let $m \asymp n^{r/(2+r)}$ and
let $k \geq n$.
Then we have
$\mathbb{E}\rho(F,\hat{F}_Z) = O(n^{-2/(2+r)})$.\\
(2) Let $\rho$ be the KS distance. Let $m\asymp n^{r/(2+r)}$.
Then
$\mathbb{E}\rho(F,\hat{F}_Z) = O\left( \min \left(\frac{\log n}{n^{2/(2+r)}},
 \sqrt{\frac{\log n}{n}} \right) \right)$.
\end{theorem}

Hence, this method achieves the minimax rate of convergence
in $L_2$ while the first data release method does not. 
This suggests that the perturbation method is preferable for 
the $L_2$ distance. 
The perturbation method does not achieve the minimax rate of 
convergence in KS distance; in fact, the
exponential mechanism based method achieves a better rate
as we shown in Section~\ref{sec::exponential} (Theorem~\ref{thm::ks}). 
We examine this method numerically in Section \ref{sec::simulations}.

Another approach to histograms is given by
\cite{MKAGV08}.
They put a 
Dirichlet $(a_1,\ldots, a_m)$ prior on the 
cell probabilities
$p_1, \ldots, p_m$
where $p_j = \mathbb{P}(X_i\in B_j)$.
The corresponding posterior is
Dirichlet $(a_1+C_1,\ldots, a_m+C_m)$.
Next they draw
$q = (q_1,\ldots, q_m)$ 
from the posterior and finally they
sample new cell counts
$D=(D_1, \ldots, D_m)$ 
from a Multinomial $(k,q)$.
Thus, the distribution of $D$ given $X$ is
$$
\mathbb{P}(D=d|X) =
\frac{\prod_{j=1}^m \Gamma(d_j + a_j + C_j)}
{\Gamma(k + n + \sum_j a_j)}.
$$

They show that differential privacy requires
$a_j + C_j \geq k/(e^{\alpha} -1)$ for all $j$.
If we take
$a_1 = a_2 = \cdots = a_m$
then this is similar to the first histogram-based data release method
we discussed in this section.
They also suggest a weakened version of
differential privacy.

\section{Exponential Mechanism}
\label{sec::exponential}

In this section we will consider the 
exponential mechanism in some detail.
We'll derive some general results about
accuracy and apply the method to the mean,
and to density estimation.
Specifically, we will show the following for exponential mechanisms:
\begin{enumerate}
\item Choosing the size $k$ of the released database is delicate.
Taking $k$ too large compromises privacy.
Taking $k$ too small compromises accuracy.
\item The accuracy of the exponential scheme
can be bounded by a simple formula.
This formula has a term that
measures how likely it is for a distribution
based on sample size $k$, to be in a small ball
around the true distribution.
In probability theory, this is known as a small ball probability.
\item The formula can be applied to several examples
such as the KS distance, 
the mean, and nonparametric density estimation
using orthogonal series.
In each case we can use our results to choose $k$
and to find the rate of convergence of an estimator
based on the sanitized data.
\end{enumerate}
In light of Theorem
\ref{thm::nogo}, we know that some assumptions are needed on $P$.
We shall assume throughout this section that
$P$ has a bounded density $p$; note that this is a weaker condition
than \eqref{eq::lipschitz}. 

Recall the exponential mechanism.
We draw the vector 
$Z=(Z_1, \ldots, Z_k)$ from $h(z|x)$
where
\begin{eqnarray}\label{eq::gxz}
h(z|x) & = & \frac{g_x(z)}{ \int_{[0,1]^k} g_x(s) ds}, \; \; \; \text{ where }
g_x(z) \; = \;  
\exp\left(- \frac{\alpha\  \rho(\hat{F}_x,\hat{F}_z)}{2 \Delta_{n,k}}\right)
\; \; \text { and }  \\ \nonumber
\Delta & \equiv & \Delta_{n,k} \; = \;
\sup_{\stackrel{x,y\in {\cal X}^n}{\delta(x,y)=1}} 
\sup_{z\in {\cal X}^k}
|\rho(\hat{F}_x,\hat{F}_z) - \rho(\hat{F}_y,\hat{F}_z)|.
\end{eqnarray}

\begin{lemma}\label{prop::delta}
For KS distance
$\Delta_{n, k} \leq \inv{n}$.
\end{lemma}

This framework is used in Blum et al. (2008).
For the rest of this section,
assume that
$Z=(Z_1,\ldots, Z_k)$
are drawn from an exponential mechanism
$Q_n$.

\begin{definition}
Let $F$ denote the cumulative distribution function on $\mathcal{X}$ 
corresponding to $P$.
Let $\hat{G}$ denote the empirical cdf
from a sample of size $k$ from $P$, 
and let
$$
R(k,\epsilon) = P^k(\rho(F,\hat{G}) \leq \epsilon).
$$
$R(k,\epsilon)$ is called the small ball probability associated with $\rho$.
\end{definition}

The following theorem bounds the accuracy of the
estimator from the sanitized data by a simple formula
involving the small ball probability.

\begin{theorem}\label{thm::exponential}
Assume that $P$ has a bounded density $p$, 
and that there exists $\epsilon_n \to 0$ such that
\begin{equation}
\label{eq::emp-large-dev}
\mathbb{P}\left(\rho(F,\hat{F}_X) > \frac{\epsilon_n}{16}\right)=
 O\left(\inv{n^c}\right)
\end{equation}
for some $c>1$.
Further suppose that $\rho$ satisfies the triangle inequality.
Let $Z=(Z_1,\ldots, Z_k)$ be drawn from $g_x(z)$ given in \eqref{eq::gxz}.
Then,
\begin{equation}\label{eq::aaaa}
\mathbb{P}\left( \rho(F,\hat{F}_Z)> \epsilon_n\right) \leq
\frac{(\sup_x p(x))^k \exp{\left(\frac{-3\alpha \epsilon_n}{16\Delta}\right)}}
{R(k,\epsilon_n/2)} +  O\left(\inv{n^c}\right).
\end{equation}
\end{theorem}

Thus, if we can choose $k=k_n$ in such a way that
the right hand side of
(\ref{eq::aaaa}) goes to 0, then the mechanism is consistent.
We now show some examples that satisfy these conditions
and we show how to choose $k_n$.

\subsection{The KS Distance}

\begin{theorem}\label{thm::ks}
Suppose that $P$ has a bounded density $p$ and 
let $B := \log \sup_x p(x) > 0$.
Let $Z=(Z_1,\ldots, Z_k)$ be drawn from $g_x(z)$ given in (\ref{eq::gxz})
with $\rho$ being the KS distance. By requiring that 
$k_n \asymp \left(\frac{3 \alpha}{B}\right)^{2/3} n^{2/3}$, 
we have for $\epsilon_n = 2 \left(\frac{B}{3\alpha}\right)^{1/3} n^{-1/3}$,
and for $\rho$ being the KS distance,
\begin{equation}
\label{eq::ks-bound-proof}
\rho(F,\hat{F}_Z) = O_P\left(\epsilon_n \right).
\end{equation}
\end{theorem}

Note that
$\rho(F,\hat{F}_Z)$ converges to 0
at a slower rate than $\rho(F,\hat{F}_X)$.
We thus see that the rate after sanitization is $n^{-1/3}$
which is slower than the optimal rate of $n^{-1/2}$.
It is an open question whether this
rate can be improved.

\subsection{The Mean}

It is interesting to consider what happens when
$\rho(F,\hat{F}_Z) = ||\mu - \overline{Z}||^2$
where
$\mu = \int x dP(x)$ and
$\overline{Z}$ is the sample mean of
$Z$.
In this case
$\Delta \leq r/n$.
Thus,
$h(u|x) \approx e^{-n ||\overline{X}-\overline{Z}||^2/(2\alpha)}$
so, approximately,
$Z_1, \ldots, Z_k \sim N(\overline{X},k\alpha/n)$.
Indeed, it suffices to take $k=1$ in this case
since then
$\overline{Z} = \overline{X} + O_P(1/\sqrt{n})$.
Thus $\overline{Z}$ converges at the same rate as $\overline{X}$.
This is not surprising: preserving
a single piece of information requires
a database of size $k=1$.

\section{Orthogonal Series Density Estimation}
\label{sec::density}

In this section, we develop an exponential scheme based on 
density estimation and we compare it to the perturbation approach.
For simplicity we take $r=1$.
Let $\{1,\psi_1,\psi_2,\ldots,\}$
be an orthonormal basis
for $L_2(0,1) = \{f:\ \int_0^1 f^2(x) dx < \infty\}$
and assume that $p\in L_2(0,1)$.
Hence
\begin{eqnarray*}
p(x) & = & 1 + \sum_{j=1}^\infty \beta_j \psi_j(x)
\;  \text{ where } \;
\beta_j = \int_0^1 \psi_j(x) p(x) dx.
\end{eqnarray*}
We assume that the basis functions are uniformly bounded
so that
\beq
\label{eq::constant}
c_0 \equiv \sup_j \sup_x |\psi_j(x)| < \infty.
\eeq
Let ${\cal B}(\gamma,C)$ denote the Sobolev ellipsoid
$$
{\cal B}(\gamma,C) =\Biggl\{ \beta = (\beta_1,\beta_2,\ldots):\ 
\sum_{j=1}^\infty \beta_j^2 j^{2\gamma} \leq C^2 \Biggr\}
$$
where $\gamma > 1/2$.
Let
$$
{\cal P}(\gamma,C) = \Biggl\{ 
p(x) = 1 + \sum_{j=1}^\infty \beta_j \psi_j(x):\ 
\beta \in {\cal B}(\gamma,C) \Biggr\}.
$$
The minimax rate of convergence
in $L_2$ norm
for ${\cal P}(\gamma,C)$ is
$n^{-2\gamma/(2\gamma+1)}$~\citep{efro:1999}.
Thus
$$
\inf_{\hat{p}} \sup_{P\in {\cal P}(\gamma,C)} E \int (\hat{p}(x) - p(x))^2 dx \geq c_1 n^{-2\gamma/(2\gamma+1)}
$$
for some $c_1>0$.
This rate is achieved by the estimator
\beq
\label{eq::estimator-x}
\hat{p}(x) = 1+\sum_{j=1}^{m_n}\hat\beta_j \psi_j(x)
\eeq
where
$m_n = n^{1/(2\gamma+1)}$ and
$\hat\beta_j = n^{-1}\sum_{i=1}^n \psi_j(X_i).$
See \cite{efro:1999}.

For a function $u \in L_2(0,1)$, let us define
$\norm{u}_{\ell_2} =  \left( \int_0^1 |u(x)|^2  dx \right)^{1/2}$,
which is a norm on $L_2(0, 1)$.
Now consider an exponential mechanism based on
\begin{eqnarray}
\label{eq::density-ell-2-distance}
\xi(X,Z) &  = &
  \left(\int\left( \hat{p}(x) - \hat{p}^*(x)\right)^2 dx\right)^{1/2}:=  
\norm{\hat{p} - \hat{p}^*}_{\ell_2}  \hspace{1cm} \text{where} \\
\label{eq::estimator-x*}
\hat{p}^*(x) & = & 1+\sum_{j=1}^{m_k}\hat\beta_j^* \psi_j(x),
\; \text{ for} \; \; \; m_k = k^{\inv{2\gamma+1}} \; \text{ and } \;
\hat\beta_j^* =  k^{-1}\sum_{i=1}^k \psi_j(Z_i).
\end{eqnarray}

\begin{lemma}
\label{prop::ell-2-distance}
Under the above scheme we have
$\Delta \leq \frac{2c_0^2 m_n}{n}$
for $c_0$ as defined in~\eqref{eq::constant}.
Hence,
\begin{equation}\label{eq::this-is-g}
g(z|x) = 
\exp\left(- \frac{\alpha \norm{\hat{p}^* - \hat{p}}_{\ell_2}}
{\Delta} \right) \leq
\exp\left(- \frac{\alpha n\norm{\hat{p}^* - \hat{p}}_{\ell_2}}
{2c_0^2 m_n} \right) \; \text{ almost surely.}
\end{equation}
\end{lemma}

\begin{theorem}\label{thm::den-exp}
Let $Z=(Z_1,\ldots, Z_k)$ be drawn from $g_x(z)$ given in (\ref{eq::this-is-g}).
Assume that $\gamma > 1$.
If we choose $k \asymp \sqrt{n}$ then
$$
\rho^2(p,\hat{p}^*) = 
O_P\left( n^{- \frac{\gamma}{2\gamma +1}} \right).
$$
\end{theorem}

We conclude that the sanitized estimator converges
at a slower rate than the minimax rate.
Now we compare this to the perturbation approach.
Let $Z=(Z_1,\ldots, Z_k)$ be an iid sample from
$$
\hat{q}(x) = 1 + \sum_{j=1}^{m_n} (\hat\beta_j + \nu_j)\psi_j(x)
$$
where
$\nu_1, \ldots, \nu_m$ are iid draws
from a Laplace distribution with 
density 
$g(\nu) = (n \alpha/(2c_0 m))e^{- n\alpha |\nu|/(c_0 m)}$.
Thus, i the notation of \ref{lemma::useful},
$R=(\nu_1, \ldots, \nu_m)$.
It follows from Lemma \ref{lemma::useful}
that, for any $k$, this preserves differential privacy.
If $\hat{q}(x) < 0$ for any $x$
then we replace $\hat{q}$ by
$\hat{q}(x) I(\hat{q}(x)>0)/\int \hat{q}(s) I(\hat{q}(s)>0)ds$
as in~\cite{HM93}.

\begin{theorem}\label{thm::den-orth-pert}
Let $Z=(Z_1,\ldots, Z_k)$ be drawn from 
$\hat{q}$.
Assume that $\gamma > 1$.
If we choose $k \geq n$, then
$$
\rho^2(p,\hat{p}_Z) = 
O_P\left( n^{- \frac{2\gamma}{2\gamma +1}} \right)
$$
where $\hat{p}_Z$ is the orthogonal series density estimator based on $Z$.
\end{theorem}

Hence, again, the perturbation technique achieves the
minimax rate of convergence
and so appears to be superior to the exponential mechanism.
We do not know if this is because
the exponential mechanism is inherently less accurate,
or if our bounds for the exponential mechanism are not tight enough.

\section{Example}
\label{sec::simulations}

\begin{figure}
\begin{center}
\includegraphics{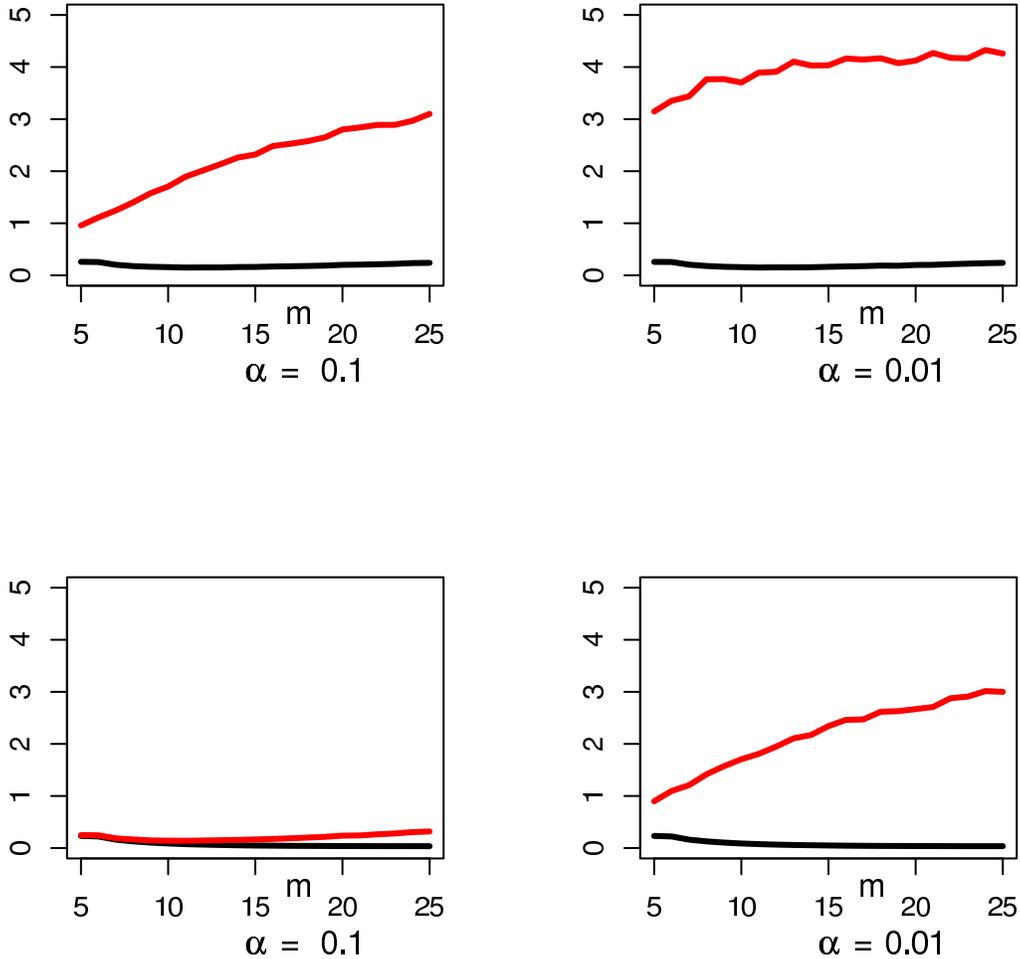}
\end{center}
\caption{Top two plots $n=100$.
Bottom two plots $n=1,000$.
Each plot shows the mean integrated squared error of the histogram.
The lower line is from the histogram based on the original data.
The upper line is based on the perturbed histogram.}
\label{fig::hist}
\end{figure}

Here we consider a small simulation study to see the effect of perturbation
on accuracy.
We focus on the histogram perturbation method with $r=1$.
We take the true density of $X$ to be a Beta(10,10) density.
We considered sample sizes $n=100$ and $n=1,000$
and privacy levels $\alpha =0.1$,
and $\alpha =0.01$.
We take $\rho$ to be squared error distance.
Figure \ref{fig::hist} shows the results of 1,000 simulations
for various numbers of bins $m$.

As expected, smaller values of $\alpha$ induce a larger
information loss which manifests itself as a larger 
mean squared error.
Despite the fact that the perturbed histogram achieves
the minimax rate, the error is substantially inflated
by the perturbation.
This means that the constants
in the risk are important, not just the rate.
Also, the risk of the sanitized histograms
is much more sensitive to the choice of the number of
cells than the original histogram is.

We repeated the simulations with a bimodal density, namely,
$p(x)$ being an equal mixture of a 
Beta(10,3) density and Beta(3,10) density.
The results turned out to be nearly identical to those above.

\section{Conclusion}
\label{sec::conclusion}

Differential privacy is an important
type of privacy guarantee when releasing data.
Our goal has been to present the idea in statistical language and then
to show that 
loss functions based on distributions and densities
can be useful for comparing privacy mechanisms.

We have seen that sampling from a histogram leads to differential privacy
as long as either the histogram is shifted away from 0 by a factor $\delta$
or if the cells are perturbed appropriately.
The latter method achieves a faster rate of convergence in $L_2$ distance.
But, the simulation showed that the risk can nonetheless be quite large.
This suggests that more work is needed to get precise
finite sample risk bounds.
Also, the choice of the smoothing parameter (number of cells in the histogram)
has a larger effect on the sanitized histogram than on the original histogram.

We also studied the exponential mechanism.
Here we derived a formula for assessing the accuracy of the method.
The formula involves small ball probabilities.
As far as we know, the connection between differential privacy
and small ball probabilities has not been observed before.

%
%
%

Minimaxity is desirable for any statistical procedure.
We have seen that in some cases the minimax rate
is achieved and in some cases it is not.
We do not yet have a complete minimax theory for differential
privacy and this is the focus of our current work.
We close with some open questions.

\begin{enumerate}
\item When is it possible for $\rho(F,\hat{F}_Z)$
to have the same rate as $\rho(F,\hat{F}_X)$?
\item When adaptive minimax methods are used,
such as adapting to $\gamma$
in Section \ref{sec::density} or when using wavelet estimation methods,
is some form of adaptivity preserved after sanitization?
\item Many statistical methods involve some sort of risk minimization.
A example is choosing a bandwidth by cross-validation.
What is the effect of sanitization on these procedures?
\item Are there other, better methods of sanitization
that preserve differential privacy?
\end{enumerate}

\section{Proofs}
\label{sec::proofs}

\subsection{Proof of Theorem~\ref{thm::power}}
\label{sec::proofofpower}
\begin{proofof2}
Without loss of generality take $i=1$.
Let $M_0(B)$ $=$ $\int Q(B|s,x_2,\ldots, x_n)$ $dP(x_2,\ldots, x_n)$ and
$M_1(B) = \int Q(B|t,x_2,\ldots, x_n) dP(x_2,\ldots, x_n)$.
By the Neyman-Pearson lemma, the highest power test
is to reject $H_0$ when
$U > u$ where
$U(z) = (d M_1/d M_0) (z)$ and
$u$ is chosen so that
$\int I(U(z) > u) dM_0(z) \leq \gamma$.
Since $(s,x_2,\ldots, x_n)$ and 
$(t,x_2,\ldots, x_n)$ differ in only one coordinate,
$M_1(B) \leq e^\alpha M_0(B)$ and so
the power is
$M_1 (U>u) \leq e^\alpha M_0(U>u) \leq \gamma e^\alpha$.
\end{proofof2}

\subsection{Proof of Lemma \ref{lemma::useful}}

For the first part simply note that
$\mathbb{P}( h(T(X,R))\in B|X=x) = 
\mathbb{P}( T(X,R)\in h^{-1}(B)|X=x) \leq
e^\alpha \mathbb{P}( T(X,R)\in h^{-1}(B)|X=x') =
e^\alpha \mathbb{P}( h(T(X,R))\in B|X=x').$

For the second part, let
$Z =(Z_1, \ldots, Z_k)$ and note that
$Z$ is independent of $X$ given $T(X,R)$.
Let $H$ be the distribution of $T(X,R)$.
Hence,
\begin{eqnarray*}
\mathbb{P}(Z\in B|X=x) &=&
\int \mathbb{P}(Z\in B|X=x,T=t) dH(t|X=x) dt\\
&=& \int \mathbb{P}(Z\in B|T=t) dH(t|X=x) dt\\
&=&  \int \mathbb{P}(Z\in B|T=t) \frac{dH(t|X=x)}{dH(t|X=x')}dH(t|X=x')\\
& \leq & e^\alpha \int \mathbb{P}(Z\in B|T=t) dH(t|X=x')\\
&=& e^\alpha \mathbb{P}(Z\in B|X=x').
\end{eqnarray*}

\subsection{Proof of Theorem \ref{thm::nogo}}
\label{sec::proofofnogo}
\begin{proofof2}
Our proof is adapted
from an argument
given in Theorem 5.1. of \cite{BLR08}.
Let $r=1$ so that ${\cal X}=[0,1]$.
Let $P = \delta_0$ where $\delta_0$ denotes a point mass at 0.
Then 
$P^n(X=X_{(0)})=1$ where
$X_{(0)} \equiv \{0,\ldots, 0\}$.
Assume that $Q_n$ is 
consistent.
Since $F(0)=1$, it follows that
for any $\delta>0$, 
$\mathbb{P}(\hat{F}_Z(0) > 1-\delta)\to 1$.
But since
$\mathbb{P}(\cdot) = \mathbb{E}_P Q_n(\cdot|X)$ and since
$P^n(X=X_{(0)})=1$,
this implies that
$Q_n(\hat{F}_Z(0) > 1-\delta|X=X_{(0)})\to 1$.

Let $v>0$ be any point in $[0,1]$ such that
$Q_n(Z=v|X=X_{(0)}) =0$.
Let
$X_{(1)} = \{v,0,\ldots, 0\}$,
$X_{(2)} = \{v,v,0,\ldots, 0\}$,
\ldots,
$X_{(n)} = \{v,v,\ldots, v\}$.
By assumption,
$Q_n(Z = X_{(j)}|X=X_{(0)}) =0$ for all $j \geq 1$.
Differential privacy implies that
$Q_n(Z = X_{(j)}|X=X_{(1)}) =0$ for all $j \geq 1$.
Applying differential privacy again implies that
$Q_n(Z = X_{(j)}|X=X_{(2)}) =0$ for all $j \geq 1$.
Continuing this way, we conclude that
$Q_n(Z = X_{(j)}|X=X_{(n)}) =0$ for all $j \geq 1$.

Next let $P=\delta_v$.
Arguing as before, we know that
$Q_n(\hat{F}_Z(v) < 1-\delta|X=X_{(n)})\to 0$.
And since $F(v-) = 0$ we also have that
$Q_n(\hat{F}_Z(v-) > \delta|X=X_{(n)})\to 0$.
Here, $F(v-) = \lim_{i\to \infty} F(v_i)$
where $v_1 < v_2 < \ldots$ and
$v_i\to v$.
Hence, for $j/n > 1-\delta$,
$Q_n(Z = X_{(j)}|X=X_{(n)}) >0$ which is a contradiction.
\end{proofof2}


\subsection{Proof of Theorem \ref{thm::hist-theorem}}
\label{sec::proofofhist}
\begin{proofof2}
Suppose that $X$ differs from $Y$ in at most one observation.
Let $\hat{f}$ denote the perturbed histogram $\hat{f}_{m, \delta}$ based on 
$X$ and let $\hat{g}_{m, \delta}$ denote the histogram based on $Y$,
such that $X$ and $Y$ differ in one entry. We also use 
$\hat{p}_j(X)$ and $\hat{p}_j(Y)$ 
for cell proportions.
Note that $|\hat{p}_j(X) - \hat{p}_j(Y)| < 1/n$ by definition.
It is clear that the maximum density ratio for a single draw $x_i$, or all $i$,
occurs in one bin $B_j$. Now consider 
$\vecx = (x_1, \ldots, x_i)$ such that for all $i =1, \ldots, k$, we have 
$x_i \in B_j \subset [0, 1]^r$ and the following bounds.

\begin{enumerate}
\item
Let $\hat{p}_j(Y) = 0$; then in order to maximize 
$\hat{f}(\vecx)/\hat{g}(\vecx)$, we let $\hat{p}_j(X) = 1/n$ and obtain
$$
\frac{\hat{f}(\vecx)}{\hat{g}(\vecx)} =
\prod_{i=1}^k \frac{\hat{f}_{m, \delta}(x_i)}
{\hat{g}_{m, \delta}(x_i)} \leq \left(\frac{ (1-\delta)m(1/n) + \delta}
{\delta} \right)^k = \left( \frac{(1-\delta)m}{n \delta} + 1 \right)^k;
$$
\item
Otherwise, we let $\hat{p}_j(Y) \geq 1/n$,
(as by definition of $\hat{p}_j$, it takes $z/n$ for non-negative 
integers $z$) and let $\hat{p}_j(X) = \hat{p}_j(Y) \pm 1/n$.
Now it is clear that in order to maximize the density ratio at $x$, we 
may need to reverse the role of $X$ and $Y$,
\begin{eqnarray*}
\max\left(\frac{\hat{g}(\vecx)}{\hat{f}(\vecx)}, \frac{\hat{f}(\vecx)}{\hat{g}(\vecx)}\right)
& \leq &
\max\left(
\left(\frac{ (1-\delta)m \hat{p}_j + \delta}
     { (1-\delta)m(\hat{p}_j - (1/n)) + \delta}\right)^k, 
\left(\frac{ (1-\delta)m(\hat{p}_j + 1/n) + \delta}
     { (1-\delta)m \hat{p}_j + \delta} \right)^k \right), \\
& \leq & \max \left(\frac{ (1-\delta) m (1/n)}
     { (1-\delta)m(\hat{p}_j - (1/n)) + \delta} + 1 \right)^k  \\
& \leq & 
\left( \frac{(1-\delta)m}{n \delta} + 1 \right)^k,
\end{eqnarray*}
where the maximum is achieved when $\hat{p}_j(Y) = 1/n$ and 
$\hat{p}_j(X) = 0$, given a fixed set of parameters $m, n, \delta$.
\end{enumerate}
Thus we have 
$$\sup_{\vecx \in ([0,1]^{r}, \ldots, [0, 1]^r)} 
\frac{\hat{f}(\vecx)}{\hat{g}(\vecx)}
\leq \left( \frac{(1-\delta)m}{n\delta} + 1\right)^k,$$
and the theorem holds.  
\end{proofof2}

\subsection{Proof of Theorem~\ref{thm::histogram}}
\label{sec::proofofhistogram}
\begin{proofof2}
Recall that $\hat{F}_Z$ denotes the empirical distribution function
corresponding to $Z = (Z_1, \ldots, Z_k)$, 
where $Z_i \in [0, 1]^r$ for all $i$ are i.i.d. 
draws from density function 
$\hat{f}_{m,\delta}(x)$ as in~\eqref{eq::reasoning} 
given $X = (X_1, \ldots, X_n)$.
Let $U$ denote the uniform cdf on $[0,1]^r$.
Given $X = (X_1, \ldots, X_n)$ drawn from a distribution whose cdf is $F$,
let $\hat{f}_m$ denote the histogram estimator on $X$ and let 
$\hat{F}_m (x) = \int_0^x \hat{f}_m(s) ds$ and
$\hat{F}_{m,\delta}(x)  = (1- \delta) \hat{F}_m(x) + \delta U(x)$.
Define $F_m(x) = \mathbb{E}(\hat{F}_m(x))$ and
$\bar{f}_m(x) = \mathbb{E}(\hat{f}_m(x))$.

The Vapnik-Chervonenkis dimension of the class of sets of the form
$\{ (-\infty,x_1] \times \cdots \times (-\infty,x_r]$ is $r$
and so by the standard Vapnik-Chervonenkis bound, we have 
for $\epsilon>0$ that
\begin{eqnarray}
\label{eq::VC-N}
\mathbb{P}\left(\sup_{t \in [0,1]^r} |\hat{F}_X(t) - F(t)| > \epsilon \right)
& \leq & 8 n^r \exp\left\{- \frac{n\epsilon^2}{32} \right\} \leq
\exp\left\{- \frac{n\epsilon^2}{64} \right\}
\end{eqnarray}
for large $n$.
Hence,
$\mathbb{E}\sup_{t \in [0,1]^r} |\hat{F}_X(t) - F(t)| =
O\left(\sqrt{\frac{r \log n}{n}} \right)$.
Given $X$, we have 
$Z_1,\ldots, Z_k\sim \hat{F}_{m,\delta}$ and so
$\mathbb{E}\sup_{[0,1]^r} |\hat{F}_Z(t) - \hat{F}_{m,\delta}(t)| =
O\left(\sqrt{\frac{r \log k}{k}} \right).$ 
Thus,
\begin{eqnarray*}
\mathbb{E}\sup_{x \in [0,1]^r} \abs{\hat{F}_Z(x) - F(x)} & \leq &
\mathbb{E}\sup_x | \hat{F}_Z(x) - \hat{F}_{m,\delta}(x)| + \mathbb{E} \sup_x | \hat{F}_{m,\delta}(x) - F(x)|\\
& \leq &
\mathbb{E}\sup_x | \hat{F}_Z(x) - \hat{F}_{m,\delta}(x)| +  \mathbb{E} \sup_x | \hat{F}_{m}(x) - F(x)| + \delta\\
& \leq &
\mathbb{E}\sup_x | \hat{F}_Z(x) - \hat{F}_{m,\delta}(x)| +  \mathbb{E} \sup_x | \hat{F}_{m}(x) - F(x)| + \delta\\
&=& O\left(\sqrt{\frac{r \log k}{k}}\right) +  \mathbb{E} \sup_x | \hat{F}_{m}(x) - F(x)| + \delta.
\end{eqnarray*}
By the  triangle inequality, we have for all $x \in [0, 1]^r$,
$$\abs{ \hat{F}_{m}(x) - F(x)}  \leq 
\abs{\hat{F}_m(x) - F_m(x)} + \abs{F_m(x) - F(x)},$$
and hence
\begin{eqnarray}
\nonumber
\mathbb{E}\sup_{x \in [0, 1]^r} \abs{\hat{F}_{m}(x) - F(x)|}
& \leq & 
\mathbb{E}\sup_{x \in [0, 1]^r} \abs{\hat{F}_m(x) - F_m(x)} + 
\mathbb{E} \sup_{x \in [0, 1]^r} \abs{F_m(x) - F(x)} \\
\label{eq::supbound}
& = & O\left(\sqrt{\frac{r \log n}{n}}\right) + 
\mathbb{E} \sup_{x \in [0, 1]^r} \abs{F_m(x) - F(x)}
\end{eqnarray}
where the last step follows from the VC bound as 
in~\eqref{eq::VC-N} for $F_m(x)$.

Next we bound $\sup_{x \in [0, 1]^r} \abs{F_m(x) - F(x)}$.
Now $F(x) = P(A)$ where
$A= \{(s_1,\ldots, s_r):\ s_i \leq x_i , i=1,\ldots , r\}$.
If $x=(j_1h,\ldots, j_rh)$ for some
integers $j_1,\ldots, j_r$ then
$F(x) - F_m(x) =0$.
For $x$ not of this form,
let $\tilde{x} = (j_1h,\ldots, j_r h)$ 
where $j_i = \lfloor  x_i/h \rfloor$.
Let
$R= \{(s_1,\ldots, s_r):\ s_i \leq \tilde{x}_i , i=1,\ldots , r\}$.
So
\begin{eqnarray}
\nonumber
F(x) - F_m(x) & = & P(A) - P_m(A) 
 = P(R) - P_m(R) + P(A \setminus R) - P_m(A \setminus R) \\
\label{eq::F-diff}
& = & P(A \setminus R) - P_m(A \setminus R)
\end{eqnarray}
where $P_m(B) = \int_B d F_m(u)$ and the set $A \setminus R$ intersects at 
most $r h /h^r$ number of cubes in $\{B_1,\ldots, B_m\}$,  given
that $\vol(A \setminus R) \leq 1 - (1-h)^r \leq r h$. 
Now by the Lipschitz condition~\eqref{eq::lipschitz}, we have
$\sup_{x \in [0, 1]^r} | p(x) - \bar{f}_m(x)| \leq  L h \sqrt{r}$ and
\begin{eqnarray}
\nonumber
\lefteqn{|P(A \setminus R) - P_m(A \setminus R)| }  \\ \nonumber
&\leq &
{\rm number\ of\ cubes\ intersecting } (A \setminus R) \times
{\rm maximum\ density\ discrepancy}\times
{\rm volume\ of\ cube}\\
\label{eq::P-diff}
& \leq &
(r h/h^r) \cdot (L h \sqrt{r}) \cdot h^r 
\leq L r^{3/2} m^{-2/r}.
\end{eqnarray}
Thus we have by~\eqref{eq::supbound},~\eqref{eq::F-diff} and~\eqref{eq::P-diff}
\begin{equation}\label{eq::asdf}
\mathbb{E}\sup_x | \hat{F}_{m}(x) - F(x)|  =
O\left(\sqrt{\frac{r \log n}{n}}\right) +  L r^{3/2}m^{-2/r}.
\end{equation}
Hence,
\begin{equation*}\label{eq::supterms}
\mathbb{E}\sup_x | \hat{F}_Z(x) - F(x)| =
O\left(\sqrt{\frac{r \log k}{k}}\right) + O\left(\sqrt{\frac{r \log n}{n}}\right) + L r^{3/2} m^{-2/r} + \delta.
\end{equation*}
Set
$m \asymp n^{r/(6+r)}$, $k \asymp m^{4/r}=n^{ 4/(6+r)}$ and $\delta = (mk/n\alpha)$ we get for all $n$ large enough,
$\mathbb{E}\sup_x | \hat{F}_Z(x) - F(x)| = O\left(\frac{\sqrt{\log n}}{n^{2/(6+r)}}\right).$
\end{proofof2}

\subsection{Proof of Theorem~\ref{thm::histogram2}}
\label{proofofthmhistgram2}
\begin{proofof2}
Let $\hat{f}_Z$ be the histogram based on $Z$ as in~\eqref{eq::lwidtsdefn}.
Then
$$
(\hat{f}_Z(u)-p(u))^2 \preceq 
(1-\delta)^2 (p(u) - \hat{f}_m(u))^2 + 
\delta^2 (p(u)-1)^2 + 
(\hat{f}_{m, \delta}(u) - \hat{f}_Z)^2 
$$
where $\preceq$ means less than, up to constants.
Hence,
$$
\mathbb{E}\int (\hat{f}_Z(u)-p(u))^2 du \preceq
 R_m + \delta^2 +
\mathbb{E}\int(\hat{f}_{m,\delta}(u) - \hat{f}_Z(u))^2 du
$$
where $R_m$ is the usual $L_2$ risk of a histogram
under the Lipschitz condition~\eqref{eq::lipschitz},
namely, $m^{-2/r} + m/n$.
Conditional on $X$,
$\hat{f}_Z$ is an unbiased estimate of $\hat{f}_m$
with integrated variance $m/k$.
So, 
$$
\mathbb{E}\int (\hat{f}_Z(u)-p(u))^2 du \preceq
m^{-2/r} + \frac{m}{n} + \delta^2 + \frac{m}{k}.
$$
Minimizing this, subject to
(\ref{eq::needed}) yields
$$
m \asymp n^{r/(2r+3)}, k\asymp n^{(r+2)/(2r+3)}, 
\delta \asymp n^{-1/(2r+3)}
$$
which yields
$
\mathbb{E}\int (\hat{f}_Z(u)-p(u))^2 du  = O(n^{-2/(2r+3)}).
$
\end{proofof2}

\subsection{Proof of Theorem \ref{thm::histogram3}}
\label{proofofthmhistgram3}
\begin{proofof2}
(1) Note that
$p - \hat{f}_Z = p - \tilde{f} + \tilde{f} - \hat{f}_Z =
p - \tilde{f} + O_P\left(\frac{m}{k}\right).$
When $k \geq n$, the latter error is lower order than the other
terms and may be ignored.
Now,
$$
p(x) - \tilde{f}(x) = 
p(x) - \hat{f}_m(x) + \hat{f}_m(x) -\tilde{f}(x).
$$
Thus
$$
\int (p(x) - \tilde{f}(x))^2 dx \preceq
\int (p(x) - \hat{f}_m(x))^2 dx + 
\int (\hat{f}_m(x) -\tilde{f}(x))^2 dx.
$$
The expected value of the first term is the usual risk,
namely, $O(m^{-2/r} + m/n)$.

For the second term, we proceed as follows.
Let
$\hat{p}_j = C_j/n$ and
$$
\hat{q}_j = \frac{ (C_j + \nu_j)_+ }{\sum_{s=1}^m (C_s + \nu_s)_+ }.
$$
We claim that
$$
\max_j |\hat{q}_j -\hat{p}_j| = O\left( \frac{ \log m }{n} \right)
$$
almost surely, for all large $n$. We have
$$
\hat{q}_j = \frac{ (C_j + \nu_j)_+ }{n}
\left(\frac{n}{\sum_{s=1}^m (C_s + \nu_s)_+ }\right) =
\frac{ (C_j + \nu_j)_+ }{n} \frac{1}{R_n}
$$
where
$R_n =(\sum_{s=1}^m (C_s + \nu_s)_+)/n$.
Now
$$
\hat{p}_j - \frac{|\nu_j|}{n} \leq
\hat{p}_j + \frac{\nu_j}{n} =
\frac{ (C_j + \nu_j) }{n} \leq
\frac{ (C_j + \nu_j)_+ }{n}
\leq
\hat{p}_j + \frac{|\nu_j|}{n}.
$$
Therefore,
$$
\left| \frac{ (C_j + \nu_j)_+ }{n} - \hat{p}_j \right| \leq
\frac{ |\nu_j|}{n} \leq \frac{M}{n}
$$
where
$M = \max \{ |\nu_1|, \ldots, |\nu_m| \}$.
Let $A>0$. 
The density for $\nu_j$ has the form
$f(\nu) = (\beta/2)e^{-\beta|\nu|}$.
So,
$$
\mathbb{P}( M > A \log m )  \leq  m \mathbb{P}( |\nu_j| > A \log m )
=  \beta m \int_{A\log m}^\infty e^{-\beta |\nu|} d\nu =
\frac{1}{m^{A \beta-1}}.
$$
By choosing $A$ large enough we have that
$M < A \log m$ a.s. for large $n$,
by the Borel-Cantelli lemma.
Therefore,
$$
\left| \frac{ (C_j + \nu_j)_+ }{n} - \hat{p}_j \right| \leq
 \frac{\log m}{n}
$$

Now we bound $R_n$.
We have
$$
1 - \frac{\sum_s  |\nu_s|}{n} \leq
1 + \frac{\sum_s  \nu_s}{n} \leq
R_n = \frac{\sum_{s=1}^m (C_s + \nu_s)_+}{n}\leq
1 + \frac{\sum_s | \nu_s|}{n}
$$
so that
$$
|R_n - 1| \leq \frac{\sum_s | \nu_s|}{n} 
\leq \frac{Mm}{n} = O\left(\frac{m\log m}{n}\right)\ \ \ a.s.
$$
Therefore,
$1/R_n = (1+O(m\log m/n))$ and thus
\begin{eqnarray*}
\hat{q}_j &=& \left(\hat{p}_j + O\left(\frac{\log m}{n}\right)\right)\   
\left( 1 + O\left( \frac{m\log m}{n}\right)\right)\\
&=&
\hat{p}_j + 
\hat{p}_j \ O\left(\frac{m\log m }{n}\right) + 
O\left(\frac{\log m}{n}\right) + O\left( \frac{m (\log m)^2}{n^2}\right).
\end{eqnarray*}

Next we claim that
$\hat{p}_j = O(1/m)$ a.s.
To see this, note that $p_j \leq C/m$, 
by definition of $C$: $1 \leq C = \sup_x p(x) < \infty$.
Hence, by Bernstein's inequality,
\begin{eqnarray*}
\mathbb{P}\left(\hat{p}_j > \frac{2C}{m}\right) &=&
\mathbb{P}\left(\hat{p}_j -p_j > \frac{2C}{m}-p_j\right) \leq
\exp\left\{-\frac{1}{2}\frac{ n ( (2C/m) - p_j )^2}{ p_j + \frac{1}{3}( (2C/m) - p_j)} \right\}\\
& \leq &
\exp\left\{-\frac{1}{2}\frac{ n C^2/m^2}{(4C/3m)} \right\} =
e^{- 3 n C/(8m)} \leq \frac{1}{n^2}
\end{eqnarray*}
for all $n \geq 16 m \log n/3C$; 
Thus $\hat{p}_j = O(1/m)$ a.s. for all large $n$.
Thus,
$\hat{q}_j - \hat{p}_j = O(\log m/n)$
almost surely for all large $n$. Hence,
$$
\mathbb{E}\int (\hat{f}_m(x) -\tilde{f}(x))^2 dx =
O\left( \frac{m\log m}{n}\right)^2.
$$
So the risk is
$$
O\left( m^{-2/r} + \frac{m}{n} + \left(\frac{m\log m }{n}\right)^2  \right) = 
O\left( m^{-2/r} + \frac{m}{n} \right),
$$
for $n \geq m \log^2 m$. This is the usual risk.
Hence,
we can choose $m \asymp n^{r/(2+r)}$ to achieve risk
$n^{-2/(2+r)}$ for all $n$ large enough.

(2) Let $\hat{F}_m$ be the cdf based on the original histogram and let
$\tilde{F}_m$ be the cdf based on the perturbed histogram.
We have
\begin{eqnarray*}
\mathbb{E}\sup_x |F(x) - \hat{F}_Z(x)| & \leq &
\mathbb{E}\sup_x |F(x) - \hat{F}_m(x)|  + \mathbb{E}\sup_x |\hat{F}_m(x) - \tilde{F}_m(x)| +
\mathbb{E}\sup_x |\tilde{F}_m(x) - \hat{F}_Z(x)|\\
& \leq & \mathbb{E}\sup_x |F(x) - \hat{F}_m(x)|  + 
\mathbb{E}\sup_x |\hat{F}_m(x) - \tilde{F}_m(x)| + 
O\left(\sqrt{\frac{r \log k}{k}}\right).
\end{eqnarray*}
Since we may take $k$ as large as we like, we can make the last 
term arbitrarily small.
From \eqref{eq::asdf},
$$
\mathbb{E}\sup_x |F(x) - \hat{F}_m(x)|   =  
O\left(\sqrt{\frac{r \log n}{n}}\right) + L r^{3/2}m^{-2/r}.
$$

%
Let $\hat{f}(x) = h^{-r}\sum_{j=1}^m \hat{p}_jI(x\in B_j)$ and 
Let $\tilde{f}(x) = h^{-r}\sum_{j=1}^m \hat{q}_jI(x\in B_j)$.
Let $x' = (u_1 h, \ldots, u_r h)$ where 
$u_i = \lceil  x_i /h \rceil, \forall i = 1, \ldots, r$.
Recall that $B_1, \ldots, B_m$ are the $m$ bins of $\mathcal{X}$ with sides
of length of $h$. 
Let $B_x$ denote the cube with the left-most corner being $0$ and 
the right-most corner being $x$.
Then for all $x$, we have
\begin{eqnarray*}
\abs{\hat{F}_m(x) - \tilde{F}_m(x)} & = & 
\abs{\int_0^x \hat{f}(s) - \tilde{f}(s) ds}  \leq 
\int_0^x \abs{\hat{f}(s) - \tilde{f}(s)} ds \\
& \leq & 
\int_0^{x'} \abs{\hat{f}(s) - \tilde{f}(s)} ds \\
& = & 
\sum_{\ell: B_{\ell} \subseteq B_{x'}}
|\hat{p}_{\ell} - \hat{q}_{\ell}| \leq
\sum_{\ell = 1}^m |\hat{p}_{\ell} - \hat{q}_{\ell}|
\end{eqnarray*}
where we use the fact that there are at most $m$ cubes.
Hence,
\begin{eqnarray*}
\mathbb{E}\sup_{x \in [0, 1]^r} |\hat{F}_m(x) - \tilde{F}_m(x)| & \leq &
\frac{m \log m}{n}
\end{eqnarray*}
where we use the fact that
$\max_j |\hat{p}_j - \hat{q}_j| = O(\log m /n)$ a.s.
So,
$$
\mathbb{E}\sup_x |F(x) - \hat{F}_Z(x)|  =
 O\left(\sqrt{\frac{r \log n}{n}}\right) + L r^{3/2}m^{-2/r}
+ O\left( \frac{m \log m}{n}\right).
$$
Setting
$m \asymp n^{r/(2+r)}$ yields
$$
\mathbb{E}\sup_x |F(x) - \hat{F}_Z(x)|  =
O\left( \min \left(\frac{\log n}{n^{2/(2+r)}},
 \sqrt{\frac{\log n}{n}} \right) \right)
$$
Hence for $r = 1$, the rate is $O\left(\sqrt{\frac{\log n}{n}}\right)$.
For $r  \geq 2$, the rate is dominated by the first term inside $O()$, 
and hence the rate is 
$O\left(\log n \times n^{-2/(2+r)}\right)$.
\end{proofof2}

\subsection{Proof of Theorem~\ref{thm::exponential}}
\label{sec::proofofexponential}
\begin{proofof2}
Let 
$B_\epsilon = \Bigl\{ u =(u_1,\ldots, u_k):\ 
\rho(F,\hat{F}_u) \leq \epsilon \Bigr\}$
where $\hat{F}_u$ is the empirical distribution
based on $u=(u_1,\ldots, u_k)\in {\cal X}^k$.
Also, let $A_n = \{ \rho(\hat{F}_X,F) \leq \epsilon_n/16 \}$.
For notational simplicity set
$\Delta = \Delta_{n, k}$. Then
\begin{eqnarray}
\nonumber
\mathbb{P}\left( \rho(F,\hat{F}_Z) > \epsilon_n\right) &=&
\mathbb{P}\left( \rho(F,\hat{F}_Z) > \epsilon_n,A_n\right) +
\mathbb{P}\left( \rho(F,\hat{F}_Z) > \epsilon_n,A_n^c\right)\\
& \leq & \nonumber
\mathbb{P}\left( \rho(F,\hat{F}_Z) > \epsilon_n,A_n\right) +
\mathbb{P}\left( A_n^c\right)\\ 
&=&
\label{eq::bnd}
\mathbb{P}\left( \rho(F,\hat{F}_Z) > \epsilon_n,A_n\right) +
O\left(\inv{n^c}\right).
\end{eqnarray}
By the triangle inequality
$\rho(\hat{F}_u,\hat{F}_X) \geq \rho(\hat{F}_u,F) - \rho(\hat{F}_X,F)$. 
Then,
\begin{eqnarray*}
\int_{B_\epsilon^c} g_x(u) du 
& = & 
\int_{B_\epsilon^c} \exp{\left(\frac{-\alpha \rho(\hat{F}_X,\hat{F}_u)}{2\Delta} 
\right)} du  \\
& \leq &
\int_{B_\epsilon^c} \exp{\left(\frac{-\alpha (\rho(\hat{F}_u,F) - \rho(\hat{F}_X,F))}
{2\Delta} \right) }du\\
&=&
\exp{\left(\frac{\alpha \rho(\hat{F}_X,F)}{2\Delta}\right)}
\int_{B_\epsilon^c} \exp{\left(\frac{-\alpha \rho(\hat{F}_u,F)}{2\Delta} \right)} du\\
& \leq &
\exp{\left(\frac{\alpha \rho(\hat{F}_X,F)}{2\Delta}\right)}
\exp{\left(\frac{-\alpha \epsilon}{2\Delta} \right)}
\int_{B_\epsilon^c}  du \\
& \leq &
\exp{\left(\frac{\alpha \rho(\hat{F}_X,F)}{2\Delta}\right)}
\exp{\left(\frac{-\alpha \epsilon}{2\Delta} \right)}.
\end{eqnarray*}
By the triangle inequality, we also have 
$
\rho(\hat{F}_u,\hat{F}_X) \leq
\rho(\hat{F}_u,F) + \rho(\hat{F}_X,F)$
and
\begin{eqnarray*}
\int g_x(u)du & \geq & \int_{B_{\epsilon/2}} g_x(u) du =
\int_{B_{\epsilon/2}}  
\exp{\left(\frac{-\alpha \rho(\hat{F}_X,\hat{F}_u)}{2\Delta}\right)} du\\
& \geq &
\exp{\left(\frac{- \alpha \rho(\hat{F}_X,F)}{2\Delta}\right)}
\int_{B_{\epsilon/2}}  
\exp{\left(\frac{-\alpha \rho(F,\hat{F}_u)}{2\Delta}\right)} du\\
& \geq &
\exp{\left(\frac{- \alpha \rho(\hat{F}_X,F)}{2\Delta}\right)}
\exp{\left(\frac{-\alpha \epsilon}{4\Delta}\right)}
\int_{B_{\epsilon/2}}  du\\
&=&
\exp{\left(\frac{-2 \alpha \rho(\hat{F}_X,F) - \alpha \epsilon}{4\Delta}\right)} 
\int_{B_{\epsilon/2}}  \frac{p(u_1)\cdots p(u_k)}{p(u_1)\cdots p(u_k)} du\\
&\geq&
\frac{
\exp{\left(\frac{-2 \alpha \rho(\hat{F}_X,F) - \alpha \epsilon}{4\Delta}\right)}}
 {(\sup_x p(x))^k}
\mathbb{P}\left( \rho(F,\hat{G}) \leq \epsilon/2\right)
\end{eqnarray*}
where $\hat{G}$ is the empirical cdf from a sample of size $k$ drawn from $P$.
Thus we have 
\begin{eqnarray*}
\label{eq::target}
\int_{B_\epsilon^c} h(u|x) du
& \leq & 
\frac{ {(\sup_x p(x))^k} \exp{\left(\frac{\alpha \rho(\hat{F}_X,F)}{\Delta}\right)}
\exp{\left(\frac{- \alpha \epsilon}{4\Delta}\right)}}
{\mathbb{P}\left( \rho(F,\hat{G}) \leq \epsilon/2\right)}.
\end{eqnarray*}
Thus, from (\ref{eq::bnd}),
\begin{eqnarray*}
\label{eq::target2}
\mathbb{P}\left( \rho(F,\hat{F}_Z) > \epsilon\right)  
& \leq &
\mathbb{P}\left(\rho(\hat{F}_X,F) \geq \frac{\epsilon}{16}\right) +
\frac{ {(\sup_x p(x))^k}
\exp{\left(\frac{- 3 \alpha \epsilon}{16\Delta}\right)}}
{\mathbb{P}\left( \rho(F,\hat{G}) \leq \epsilon/2\right)} \\
& = &
\frac{ {(\sup_x p(x))^k}
\exp{\left(\frac{- 3 \alpha \epsilon}{16\Delta}\right)}}
{\mathbb{P}\left( \rho(F,\hat{G}) \leq \epsilon/2\right)}
+  O\left(\inv{n^c}\right).
\end{eqnarray*}
Thus the theorem holds.
\end{proofof2}

\subsection{Proof of Lemma~\ref{prop::delta}}
\label{sec::proofoflemmadelta}
\begin{proofof2}
{\bf Proof of Lemma \ref{prop::delta}.}
We start with KS,
By the triangle inequality, we have 
for all $z \in \X^k$ and for all
$x,y \in \X^n$,
\ben
\abs{\rho(\hat{F}_x,\hat{F}_z) - \rho(\hat{F}_y,\hat{F}_z)} 
& \leq &
\rho(\hat{F}_x, \hat{F}_y).
\een
Notice that changing one entry in $x$ will change
$\hat{F}_x (t)$ by at most $\inv{n}$ at any $t$ by definition,
that is,
\ben
\sup_{t \in [0, 1]^r} |\hat{F}_x(t) - \hat{F}_y(t)| & = & \inv{n}.
\een
Thus the conclusion holds for the KS-distance.
\end{proofof2}

\subsection{Proof of Theorem~\ref{thm::ks}}
\label{sec::proofofcvm}

We need the following small ball result; see \cite{lishao:2001}.
\begin{theorem}
Let $r \geq 3$, and $\{X_t, t \in [0, 1]^r\}$ be the Brownian sheet. 
Then there exists $0< C_r < \infty$ such that 
for all $0< \epsilon \leq 1$,
\begin{eqnarray*}
\log \mathbb{P} \left(\sup_{t \in [0, 1]^r} |X_t|\leq \epsilon\right) 
\geq - C_r \epsilon^{-2} \log^{2r-1}(1/\epsilon)
\end{eqnarray*}
where
$C_r$ depends only on $r$. The same bound holds for a Brownian bridge.
\end{theorem}

\begin{proofof}{\textnormal{theorem~\ref{thm::ks}}}
The Vapnik-Chervonenkis dimension of the class of sets of the form
$\{ (-\infty,x_1] \times \cdots \times (-\infty,x_r]$ is $r$
and so by the standard Vapnik-Chervonenkis bound, we have 
for $\epsilon_n, k_n$ as specified in the theorem statement,
\begin{eqnarray}
\nonumber
\mathbb{P}\left(\sup_{[0,1]^r} |\hat{F}_X(t) - F(t)| > \frac{\epsilon_n}{16} \right)
& \leq & 8 n^r \exp\left\{- \frac{n (\epsilon_n/16)^2}{32} \right\} \\ 
\nonumber
& \leq & 
8\exp\left\{-c_5 \left(\frac{B}{3\alpha}\right)^{2/3} n^{1/3}
+ r \log n \right\} \\ \nonumber
& = & 
8\exp\left\{- c_6 \sqrt{k_n} \left(\frac{B}{3\alpha}\right) + c_7 r \log k_n\right\} \\
\label{eq::VC-bound}
& = & 
8\exp\left\{- C_2 \sqrt{k_n} \left(\frac{B}{3\alpha}\right)\right\}
\end{eqnarray} 
for some constants $c_5, c_6, c_7, C_2 > 0$ for $n$ large enough.
Thus (\ref{eq::emp-large-dev}) holds.
Now we compute the small ball probability.
Note that $\sqrt{k}(\hat{F}_k -F)$ converges to a Brownian bridge $B_k$
on $[0,1]^r$. 
More precisely, from~\cite{CR75}
there exist a sequence of Brownian bridges $B_k$ such that
\beq
\label{eq::strong-approx}
\sup_t |\sqrt{k}(\hat{F}_k -F)(t) - B_k(t)| = 
O\left( \frac{ (\log k)^{3/2}}{k^\gamma}\right)\ \ \ \ \text{ a.s. }
\eeq
where $\gamma = 1/(2(r+1))$. It is clear that the RHS of 
\eqref{eq::strong-approx} is $o(1)$ a.s. given a fixed $r$.
Hence we have for $k = k_n$ and
$\epsilon_n$ as chosen in the theorem statement, and 
for all $\epsilon \geq \epsilon_n$, it holds that 
\begin{eqnarray}
\nonumber
\log \mathbb{P} (\sup_t |\hat{F}_Z(t) - F(t)| \leq \epsilon/2) 
&=&
\log \mathbb{P} (\sup_t \sqrt{k}|\hat{F}_Z(t) - F(t)| 
\leq \sqrt{k} \epsilon/2)\\
& \geq &
\label{eq::bound1}
\log \mathbb{P} \left(\sup_t |B_k(t)| \leq \sqrt{k} \epsilon - 
O\left(k^{-\gamma} (\log k)^{3/2} \right)\right) \\
& \geq &
\label{eq::bound2}
\log \mathbb{P} \left(\sup_t |B_k(t)| \leq \frac{\sqrt{k} \epsilon}{4}
\right)
\end{eqnarray}
for all large $n$,
where \eqref{eq::bound1} follows from~\eqref{eq::strong-approx}
and~\eqref{eq::bound2} holds given that 
$\sqrt{k} \epsilon \geq \sqrt{k_n} \epsilon_n \geq c$ for some constant
$c > 1/2$ due to our choice of $k_n$ and $\epsilon_n$.
Also, $\Delta \leq 1/n$ for KS distance. 
Hence, by Theorem~\ref{thm::exponential} and~\eqref{eq::VC-bound}, we have
for $B = \log \sup_x p(x) > 0$, 
\begin{eqnarray}
\nonumber
\lefteqn{\mathbb{P}\left( \rho(F,\hat{F}_Z) > \epsilon_n \right) }\\ \nonumber
& \leq & 
C_0
\exp\left\{ -n \left(\frac{3 \alpha \epsilon_n}{16} -
\frac{B k_n}{n} - 
\frac{C_1 |\log (\sqrt{k_n} \epsilon_n/4)|^{2r-1}}{n k_n \epsilon_n^2}
 \right)\right\} + 
8 \exp\left\{- C_2 \frac{B \sqrt{k_n}}{3\alpha} \right\} \\
& \leq & 
\label{eq::dominating}
C_0 \exp(- C_3 B k_n /2) + 
8 \exp\left\{ - C_2\left(\frac{B}{3\alpha}\right) \sqrt{k_n} \right\} 
\to 0
\end{eqnarray}
for some constants $C_0, C_1, C_2$ and $C_3$,
where \eqref{eq::dominating} holds when we take w.l.o.g. 
$k_n = \inv{16} \left(\frac{3 \alpha}{B}\right)^{2/3} n^{2/3}$ and 
$\epsilon_n \geq 2 \left(\frac{B}{3 \alpha}\right)^{1/3} n^{-1/3}$, given that
$\epsilon_n \geq  2 \left(\frac{B}{3 \alpha}\right)^{1/3} n^{-1/3}
= \frac{32 k_n B}{3 n \alpha}$ and hence 
$\frac{3 \alpha \epsilon_n}{16} \geq  \frac{2B k_n}{n}$.
Thus the result follows.
\end{proofof}

\begin{remark}
The constants taken in the proof are arbitrary; indeed,
when we take $k_n = C_4 \left(\frac{3 \alpha}{B}\right)^{2/3} n^{2/3}$ and 
$\epsilon_n = 32 C_4 \left(\frac{B}{3 \alpha}\right)^{1/3} n^{-1/3}$ with
some constant $C_4 \geq 1/16$,~\eqref{eq::dominating} will hold
with slightly different constants $C_2, C_3$. 
For $k_n$ and $\epsilon_n$ as chosen above, 
it holds that $\sqrt{k_n} \epsilon_n \asymp 1$. 
\end{remark}

\subsection{Proofs for Lemma~\ref{prop::ell-2-distance} and 
Theorem~\ref{thm::den-exp}}
\label{sec::proofofdenexp}
Throughout this section, we let $\hat{p}_X$ denote the estimator 
as defined in~\eqref{eq::estimator-x},
which is based on a sample of size $n$ drawn independently from $F$; 
Similarly, we let $\hat{p}_k$ denote the same estimator 
based on an i.i.d. sample $(Y_1, \ldots, Y_k)$ of size $k$ drawn from $F$,
with $m_k = k^{1/(2 \gamma + 1)}$ replacing $m_n$ and 
$\hat\beta_j = k^{-1}\sum_{i=1}^k \psi_j(Y_i)$ in~\eqref{eq::estimator-x}.
We let $\hat{p}_Z$ denote the estimator as in~\eqref{eq::estimator-x*},
based on an i.i.d. sample $Z = (Z_1, \ldots, Z_k)$ of size $k$ drawn 
from $g_x(z)$ as in~\eqref{eq::this-is-g}.

{\bf Proof of Lemma \ref{prop::ell-2-distance}.}
\begin{proofof2}
Without loss of generality,
let $X = (x, X_2,\ldots, X_n)$ and $Y = (y, X_2, \ldots, X_n)$ 
so that $\delta(X, Y) = 1$ and let $Z \in \X^k$.
Recall that 
\ben
\xi(X,Z) & = & 
\left(\int\left( \hat{p}_X(x) - \hat{p}_Z(x) \right)^2 dx\right)^{1/2}, \\
\xi(Y,Z) & = & 
\left(\int\left( \hat{p}_Y(x) - \hat{p}_Z(x) \right)^2 dx\right)^{1/2}.
\een
In particular, let us define $u = \hat{p}_X - \hat{p}_Z$
and $v = \hat{p}_Y - \hat{p}_Z$ and thus
\silent{
\ben
\norm{u}_{\ell_2} = \left(\int\abs{\hat{p}_X(x) - \hat{p}_Z(x)}^2 dx
\right)^{1/2},   \hspace{0.5cm} \text{and} \; \;x
\norm{v}_{\ell_2} = \left(\int\abs{\hat{p}_Y(x) - \hat{p}_Z(x)}^2 dx\right)^{1/2} \\
\een
Thus we have that }
\ben
\abs{\xi(X, Z) - \xi(Y, Z)} & = & 
\abs{\left(\int\left( \hat{p}_X(x) - \hat{p}_Z(x) \right)^2 dx\right)^{1/2} - 
\left(\int\left( \hat{p}_Y(x) - \hat{p}_Z(x) \right)^2 dx\right)^{1/2}} \\
& = & 
\label{eq::triangle-ineq}
\abs{\norm{u}_{\ell_2} - \norm{v}_{\ell_2}} \leq  
\norm{u - v}_{\ell_2} \\
& = & 
\norm{\hat{p}_X - \hat{p}_Z - (\hat{p}_Y - \hat{p}_Z)}_{\ell_2} =
\norm{\hat{p}_X - \hat{p}_Y}_{\ell_2} \leq
\frac{2c_0^2 m_n}{n},
\een
where the first inequality is due to the triangle inequality for the 
$\norm{.}_{\ell_2}$ and the last step is due to
\begin{eqnarray*}
\abs{\hat{p}_X(x) - \hat{p}_Y (x)} & = & 
\inv{n}\abs{\sum_{j=1}^{m_n}
\left(\sum_{i=1}^n \psi_j(X_i) - \sum_{i=1}^n \psi_j(Y_i) \right) \psi_j(x)} \\
&=&
\inv{n}\abs{\sum_{j=1}^{m_n}\left(\psi_j(X_1) - \psi_j(Y_1) \right) \psi_j(x)} \\
& \leq & 
\frac{1}{n}\sum_{j=1}^{m_n} (|\psi_j(X_1)| + |\psi_j(Y_1)|) |\psi_j(x)| \leq
\frac{2c_0^2 m_n}{n}.
\end{eqnarray*}
Hence $\Delta \leq \frac{2c_0^2 m_n}{n}$.
\end{proofof2}

{\bf Proof of Theorem \ref{thm::den-exp}.}
\begin{proofof2}
For $u = (u_1,\ldots, u_k) \in {\cal X}^k$, we let 
$$\hat{p}_u(x) = 1+\sum_{j=1}^{m_k}\hat\beta_j \psi_j(x),$$
\hspace{1cm}
where $m_k = k^{\inv{2\gamma+1}}$ and 
$\hat\beta_j = k^{-1}\sum_{i=1}^k \psi_j(u_i)$.

Let $\hat{F}_u$ be the empirical distribution based on $u$.
Our proof follows that of Theorem~\ref{thm::exponential},
with
$$
\rho(F, \hat{F}_u) = \norm{p - \hat{p}_u}_{\ell_2}
\text { and } \; \; 
\rho(F_X, \hat{F}_u) = \norm{\hat{p}_X - \hat{p}_u}_{\ell_2}
$$
as defined in~\eqref{eq::density-ell-2-distance} for $X = (X_1, \ldots, X_n)$.
Now
$$B_\epsilon = \Bigl\{ u =(u_1,\ldots, u_k):\ 
\norm{p - \hat{p}_u}_{\ell_2} < \epsilon \Bigr\}.$$
Thus the corresponding triangle inequalities that we use  to
replace that in Theorem~\ref{thm::exponential} are:
\ben
\norm{\hat{p}_u - \hat{p}_X}_{\ell_2} & \geq &
\norm{\hat{p}_u - p}_{\ell_2} - \norm{\hat{p}_X - p}_{\ell_2} 
\text { and } \\
\norm{\hat{p}_u - \hat{p}_X}_{\ell_2} & \leq &
\norm{\hat{p}_u - p}_{\ell_2} + \norm{p - \hat{p}_X}_{\ell_2}.
\een
Standard risk calculations show that
(\ref{eq::emp-large-dev}) holds for some $c>0$ with
$\rho(F, \hat{F}_X)$ being replaced with 
$\norm{\hat{p}_X - p}_{\ell_2}$.
That is, by Markov's inequality,
$$
\mathbb{P}(\norm{\hat{p}_X- p}_{\ell_2} > \epsilon) \leq 
\frac{\mathbb{E}\norm{\hat{p}_X-p}_{\ell_2}^2}{\epsilon^2}
$$
and 
(\ref{eq::emp-large-dev}) follows from the polynomial decay of the
mean squared error
$\mathbb{E}||\hat{p}_X-p||^2$.
Thus, from (\ref{eq::bnd}), for $\hat{p}_Z = \hat{p}^*$ as 
in \eqref{eq::estimator-x*},
\begin{eqnarray*}
\label{eq::target3-den}
\mathbb{P}\left( \norm{p - \hat{p}_Z}_{\ell_2}) > \epsilon\right)  
& \leq &
\mathbb{P}\left(\norm{\hat{p}_X- p}_{\ell_2}\geq \frac{\epsilon}{16}\right) +
\frac{ {(\sup_x p(x))^k}
\exp{\left(\frac{- 3 \alpha \epsilon}{16\Delta}\right)}}
{\prob{\norm{p - \hat{p}_k}_{\ell_2} \leq \epsilon/2}}\\
& = &
\frac{ {(\sup_x p(x))^k}
\exp{\left(\frac{- 3 \alpha \epsilon}{16\Delta}\right)}}
{\mathbb{P}\left( \norm{p - \hat{p}_k}_{\ell_2} \leq \epsilon/2\right)}
+  O\left(\inv{n^c}\right).
\end{eqnarray*}

We need to compute the small ball probability.
Recall that $\hat{p}_k$ denote the estimator based on a sample of size $k$.
By Parseval's relation,
$$
\int (p(x) - \hat{p}_k(x))^2 dx =
\sum_{j=1}^{m_k} (\hat\beta_j - \beta_j)^2 + \sum_{m_k+1}^\infty \beta_j^2 \leq
\sum_{j=1}^{m_k} (\hat\beta_j - \beta_j)^2 + c k^{-2\gamma/(2\gamma+1)}.
$$
Let
$U_i = (\psi_1(X_i) - \beta_1, \ldots, \psi_{m_k}(X_i) - \beta_{m_k})^T$
and
$Y_i = \Sigma_k^{-1/2} U_i$
where
$\Sigma_k$ is the covariance matrix of $U_i$.
Hence,
$Y_i$ has mean 0 and identity covariance matrix.
Let $\lambda_k$ denote the largest eigenvalue of $\Sigma_k$.
From Lemma \ref{lemma::eigen} below,
$\lambda = \limsup_{k\to\infty} \lambda_k < \infty$.
Let
$Q = \sum_{j=1}^{m_k} (\hat\beta_j - \beta_j)^2$ and let
$S = k^{-1/2}\sum_{i=1}^k Y_i$.
Then, for all large $k$, and any $\delta>0$,
\begin{eqnarray*}
\mathbb{P}(Q \leq \delta^2) &=&
\mathbb{P}(S^T \Sigma_k S \leq k\delta^2)
 \geq  \mathbb{P}\left(S^T S \leq \frac{k\delta^2}{\lambda_k}\right)
 \geq  \mathbb{P}\left(S^T S \leq \frac{k\delta^2}{2\lambda}\right).
\end{eqnarray*}
From Theorem 1.1 of \cite{bentkus:2003}
we have that
\begin{eqnarray*}
\sup_c \left|
\mathbb{P}\left(S^T S \leq c\right) -
\mathbb{P}(\chi^2_{m_k}\leq c)\right| = 
O\left(\sqrt{\frac{m_k^3}{k}}\right) = 
O\left( k^{-(\gamma-1)/(2\gamma+1)}\right).
\end{eqnarray*}
%
%
Next we use the fact 
(see \cite{rohde:2008} for example)
that
$\mathbb{P}(\chi^2_m \leq m + a) \geq
1 - e^{-a^2/(4(m+a))}$.
Let 
$k = \sqrt{n}$,
$\epsilon_n = c_1 n^{-\gamma/(2\gamma+1)}$
where
$c_1 \geq 4(2 \lambda +1) (C^2 + 1)$
$$
a = 
\frac{k(\epsilon_n/4 - C^2 k^{-2\gamma/(2\gamma+1)})}{2\lambda} - m_k
\geq (C^2 +1) n^{1/2 (2\gamma+1)} - m_k 
\geq C^2 m_k,
$$
since $m_k = k^{\inv{2\gamma+1}} = n^{1/2(2 \gamma +1)}$.
We see that for all large $k$
\begin{eqnarray*}
\mathbb{P}\left(\norm{p - \hat{p}_k}_{\ell_2} \leq \frac{\sqrt{\epsilon_n}}{2} \right) 
& = &
\mathbb{P}\left(\int (p(x) - \hat{p}_k(x))^2 dx \leq \frac{\epsilon_n}{4}\right) \\
& \geq &
\mathbb{P}\left(\sum_{j=1}^{m_k} (\hat\beta_j - \beta_j)^2 \leq 
\frac{\epsilon_n}{4} - C^2 k^{-2\gamma/(2\gamma+1)}\right) \\
&=&
\mathbb{P}\left( \chi^2_{m_k} \leq 
\frac{k(\epsilon_n/4 - C^2 k^{-2\gamma/(2\gamma+1)})}{2\lambda}\right) - 
O\left( k^{-(\gamma-1)/(2\gamma+1)}\right)\\
& \geq &
1 - \exp\left(\frac{-a^2}{4(m_k+a)}\right)  - O\left( k^{-(\gamma-1)/(2\gamma+1)}\right) \\
&\geq & \frac{1}{2}  - O\left( k^{-(\gamma-1)/(2\gamma+1)}\right).
\end{eqnarray*}
Hence
\begin{eqnarray*}
\label{eq::target2-den}
\mathbb{P}\left( \norm{p - \hat{p}_Z}_{\ell_2}) > \sqrt{\epsilon_n} \right)  
& \leq &
\mathbb{P}
\left(\norm{\hat{p}_X - p}_{\ell_2}\geq \frac{\sqrt{\epsilon_n}}{16}\right) +
\frac{ {(\sup_x p(x))^k}
\exp{\left(\frac{- 3 \alpha \sqrt{\epsilon_n}}{16\Delta}\right)}}
{\prob{\norm{p - \hat{p}_k}_{\ell_2} \leq \sqrt{\epsilon_n}/2}}\\
& = &
\frac{ {(\sup_x p(x))^k}
\exp{\left(\frac{- 3 \alpha \sqrt{\epsilon_n}}{16\Delta}\right)}}
{\mathbb{P}\left( \norm{p - \hat{p}_k}_{\ell_2} \leq \sqrt{\epsilon_n}/2\right)}
+  O\left(\inv{n^c}\right) \\
& \leq & 
\frac{ {(\sup_x p(x))^k}
\exp{\left(\frac{- 3 \alpha n \sqrt{\epsilon_n}}{32c_0^2 m_n}\right)}}
{\mathbb{P}\left( \norm{p - \hat{p}_k}_{\ell_2} \leq 
\sqrt{\epsilon_n}/2\right)} +  O\left(\inv{n^c}\right)
\end{eqnarray*}
and so for $\gamma > 1$,
\ben
\mathbb{P} (\int (\hat{p}_Z - p)^2 \leq \epsilon_n) 
& \leq & c_2 \exp\left(k \log \sup_x p(x)\right) 
\exp\left(\frac{-3 \sqrt{c_1} \alpha n}
{n^{1/(2\gamma+1)} n^{\gamma/2 (2\gamma + 1)}}\right) \\ 
& = &
c_2 \exp\left(n^{1/2} \log \sup_x p(x) -\alpha c_3
n^{\left(\frac{3 \gamma}{2(2\gamma+1)} \right)} \right) \\
& = &
c_2 \exp\left(-\alpha c_4 n^{\left(\frac{3 \gamma}{2(2\gamma+1)} \right)} \right)
\rightarrow 0,\ \ \ 
\een
as $n \rightarrow \infty$ since $\frac{3 \gamma}{2(2\gamma+1)} > 1/2$,
where $c_2, c_3, c_4$ are some constants.
Hence the theorem holds.
\end{proofof2}

\begin{lemma}
\label{lemma::eigen}
Let
$\lambda = \limsup_{k\to\infty} \lambda_k$.
Then 
$\lambda < \infty$.
\end{lemma}

\begin{proof}
Recall that the orthonormal basis is
$\psi_0,\psi_1, \ldots,$
where $\psi_0=1$ and
$\psi_j(x) = \sqrt{2}\cos (\pi j x)$.
Also
$p(x) = 1 +\sum_{j=1}^\infty \beta_j \psi_j(x)$
and $\sum_j \beta_j^2 j^{2\gamma} < \infty$.
Note that
$\sum_{j=1}^{\infty} |\beta_j|^k = O(1)$
for $k\geq 1$;
see \cite{efro:1999}.
Note that
$\Sigma_k$ is the covariance matrix of
$\hat\beta$ times $n$.
We will use the standard identities
$\cos^2 (u) = (1+\cos(2u))/2$ and
$\cos(u)\cos(v) = \frac{\cos(u-v) + \cos(u+v)}{2}.$
It follows that
$\psi_j^2(x) = 1 + \frac{1}{\sqrt{2}} \psi_{2j}(x)$
and
$\psi_j(x)\psi_k(x) = \frac{\psi_{j-k}(x) + \psi_{j+k}(x)}{\sqrt{2}}.$
Now $\mathbb{E}(\hat\beta_j) = \beta_j$.
And
$$
n {\rm Var}(\hat\beta_j) =
{\rm Var}(\psi_j(X)) = 
\int \psi_j^2(x) p(x) dx - \beta_j^2.
$$
Now
$\int \psi_j^2(x) p(x) dx =
\int \psi_j^2(x) (1 + \sum_{\ell=1}^\infty \beta_\ell \psi_\ell(x))dx
= 1 + \sum_{\ell=1}^\infty \beta_\ell \int \psi_\ell(x) \psi_j^2(x)dx
=1 + \frac{1}{2}\sum_{\ell=1}^\infty \beta_\ell \int \psi_\ell(x) 
\left(1 + \frac{\psi_{2j}(x)}{\sqrt{2}}\right)dx=
1 + \frac{\beta_{2j}}{\sqrt{2}}.$
Thus,
$\Sigma_{jj} = 1 + \frac{\beta_{2j}}{\sqrt{2}} - \beta_j^2.$
Now consider $j\neq k$.
Then
\begin{eqnarray*}
\mathbb{E}(\psi_j(X)\psi_k(X)) &=&
\int \psi_j(x)\psi_k(x) p(x) dx \\
&=& \sum_\ell \beta_\ell \int \psi_j(x)\psi_k(x) dx \\
&=&
\beta_j \int \psi_j^2(x) \psi_k(x) dx + \beta_k \int \psi_k^2(x) \psi_j(x) dx + 
\sum_{\ell\neq j,k} \beta_\ell \int \psi_j(x)\psi_k(x)  \psi_{\ell}(x) dx \\
&=&
\frac{\beta_j}{\sqrt{2}}\int \psi_{2j}(x) \psi_k(x) dx + 
\frac{\beta_k}{\sqrt{2}}\int \psi_{2k}(x) \psi_j(x) dx\\
&&\ \ \ \ +
\frac{1}{\sqrt{2}}
\sum_{\ell\neq j,k} \beta_\ell \int(\psi_{j-k}(x) + \psi_{j+k}(x))\psi_\ell (x) \\
&=&
\frac{\beta_j}{\sqrt{2}}I(2j=k) + \frac{\beta_k}{\sqrt{2}}I(2k=j) \\
&&\ \ \ +
\frac{\beta_{\ell}}{\sqrt{2}}I(\ell= |j-k| \; \& \; j \not= 2k) +
\frac{\beta_{\ell}}{\sqrt{2}}I(\ell= j+k) \\
&=&
\frac{\beta_k}{\sqrt{2}}I(2k=j) +
\frac{\beta_{|j-k|}}{\sqrt{2}}I(j \not= 2k) +
\frac{\beta_{j+k}}{\sqrt{2}} \\
&=&
\frac{\beta_{|j-k|}}{\sqrt{2}} + \frac{\beta_{j+k}}{\sqrt{2}},
\end{eqnarray*}
where we used the fact that $\psi_{-j}(x) = \psi_{j}(x)$
for all $j = 1, 2, \ldots$
and $\int \psi_j(x) dx = 0$ for all $j >0$.
So, we have for all $j \in \{1, \ldots, p\}$, 
\begin{eqnarray*} 
\sum_{k = 1}^{p} |\Sigma_{jk}|  & = & 
|\Sigma_{jj}| + \sum_{j \not= k} 
\abs{\frac{\beta_{|j-k|}}{\sqrt{2}} + \frac{\beta_{j+k}}{\sqrt{2}} 
- \beta_j\beta_k} \\
& \leq &
1 + \abs{\frac{\beta_{2j}}{\sqrt{2}}}
+ |\beta_j| \sum_{k}  |\beta_k| + \Sigma_{j \not= k} 
\abs{\frac{\beta_{|j-k|}}{\sqrt{2}}} 
+ \abs{\frac{\beta_{j+k}}{\sqrt{2}}}\\
& \leq &
1 + \abs{\frac{\beta_{2j}}{\sqrt{2}}}
+ (|\beta_j| + \sqrt{2}) \sum_{k=1}^{\infty} |\beta_k| \\
& = & O(1).
\end{eqnarray*}
Hence,
${\rm limsup}_{k\to\infty} \lambda_{\max}(\Sigma_k)
 \leq  \norm{\Sigma_k}_{\infty} = O(1)$
and the lemma holds.
\end{proof}

\subsection{Proof of Theorem \ref{thm::den-orth-pert}}
\begin{proofof2}
The proof is similar to the proof of
Theorem \ref{thm::histogram3},
so we provide a short outline.
In particular, the effect of truncation can be shown to be negligible
as in the proof of Theorem \ref{thm::histogram3}.
We have
$p-\hat{p}_Z = p- \hat{q} + \hat{q} -\hat{p}_Z=
p- \hat{q} + O_P(m/k)$
and the latter term is negligible for $k \geq n$.
Now
$p - \hat{q} = p - \hat{p} + \hat{p} - \hat{q}$.
The term $p-\hat{p}$ is the usual error term
and contributes $O(n^{-2\gamma/(2\gamma + 1)})$ to the risk.
For the second term,
$\int (\hat{p}-\hat{q})^2 = \sum_{j=1}^m \nu_j^2 = O_P(m/n) =
O_P(n^{-2\gamma/(2\gamma + 1)})$.
\end{proofof2}

\bibliography{revision}

\end{document}